\documentclass[preprint,12pt]{elsarticle}

\usepackage{tabularx}

\usepackage{amsmath,amsfonts}
\usepackage{algorithmic}
\usepackage{algorithm}
\usepackage{array}
\usepackage{float}
\usepackage{textcomp}
\usepackage{stfloats}
\usepackage{url}
\usepackage{verbatim}
\usepackage{graphicx}
\usepackage{amsmath, xparse}
\usepackage[utf8]{inputenc}
\usepackage{xcolor}
\usepackage{amsmath,amssymb,amsfonts}
\usepackage{algorithmic}
\usepackage{graphicx}
\usepackage{textcomp}
\usepackage{bm}
\usepackage{float}
\usepackage{xcolor}
\usepackage{soul}
\usepackage{xfrac}
\usepackage{subfigure}

\definecolor{LightGray}{gray}{0.8}
\definecolor{Orange}{rgb}{1.0, 0.31, 0.0}
\definecolor{Green}{rgb}{0.3, 1.0, 0.3}
\definecolor{Blue}{rgb}{0.75,0.75,1}

\newcommand{\robs}{\rho_{\rm obs}}
\newcommand{\rcyl}{\rho_{\rm cyl}}

\newcommand{\rcri}{\rho_{\rm cri}}
\newcommand{\rfil}{\rho_{\rm fil}}

\usepackage{amssymb}

\journal{}
\biboptions{sort&compress}
\begin{document}

\begin{frontmatter}



\title{Fundamentals of a\\ Null Field Method-Surface Equivalence Principle Approach for Scattering by Dielectric Cylinders}


\author[label1]{Minas Kouroublakis}

\author[label1]{Nikolaos L. Tsitsas}

\affiliation[label1]{organization={School of Informatics, Aristotle University of Thessaloniki},
            city={Thessaloniki},
            postcode={54124}, 
            country={Greece}}

\author[label3]{George Fikioris}

\affiliation[label3]{organization={School of Electrical and Computer Engineering, National Technical University of Athens},
            city={Athens},
            postcode={15780}, 
            country={Greece}}

\begin{abstract}
The null-field method (NFM) and the method of auxiliary sources (MAS) have been both used extensively for the numerical solution of boundary-value problems arising in diverse applications involving propagation and scattering of waves. It has been shown that, under certain conditions, the applicability of MAS may be restricted by issues concerning the divergence of the auxiliary currents,  manifested by the appearance of exponentially large oscillations. In this work, we combine the NFM with the surface equivalence principle (SEP) and investigate analytically the convergence properties of the combined NFM-SEP with reference to the problem of (internal or external) line-source
excitation of a dielectric cylinder. Our main purpose is to prove that (contrary to
the MAS) the discrete NFM-SEP currents, when properly normalized, always converge to the corresponding continuous current densities, and thus no divergence and oscillations phenomena appear. The theoretical analysis of the NFM-SEP is accompanied by detailed comparisons with the MAS as well as with representative numerical results illustrating the conclusions.
\end{abstract}



\begin{keyword}
Extended Boundary Condition \sep Null-Field Method \sep Extended Integral Equation \sep Surface Equivalence Principle \sep Method of Auxiliary Sources \sep Scattering \sep Dielectric Cylinder \sep Convergence \sep Oscillations



\end{keyword}

\end{frontmatter}


\section{Introduction}

The numerical solution of acoustic, electromagnetic and elastic boundary-value problems involving cylindrical geometries lies at the core of many applications, including, e.g., the analysis of scattering by periodic structures \cite{JCompPhysPeriodicCylindrical}, the modeling and design of antennas and waveguiding systems \cite{RadioScienceAntennas}, and the investigation of shielding or cloaking devices \cite{JCompPhysCloaking}. In many cases, the solution of a boundary-value problem is facilitated by the application of the so-called \emph{surface equivalence principle} (SEP) \cite{chew2022, Harrington_2001, jin2015}. 
Then, surface current densities are defined on each participating boundary, which are referred to as \emph{equivalent current densities}. Imposing the boundary conditions leads to a system of surface integral equations (SIEs), of equal number to that of the unknown current densities. 
The numerical method that is typically used to solve SIEs is the Method of Moments (MoM) \cite{Harrington_1993}. However, MoM encounters several difficulties, the main one being the evaluation of singular integrals in order to fill the diagonal elements of the pertinent impedance matrix. The singularities in the integrals stem from the fact that the integration domain coincides with the physical boundary.

The above mentioned limitations (of evaluating singular integrals in MoM) led to the development of alternative methods in which the support of the current distribution is displaced from the domain on which the boundary conditions are applied \cite{waterman1965, waterman1969, wriedt2018}. To this end, two main alternative approaches have been established: 
\begin{enumerate}
\item Define equivalent current densities on displaced surfaces; enforce boundary conditions on the physical surfaces.
\item Define equivalent current densities on the physical boundaries; enforce boundary conditions on displaced surfaces. Such boundary conditions are referred to as \emph{extended boundary conditions} (EBCs) and the resulting integral equations as \emph{extended integral equations} (EIEs).
\end{enumerate}

One of the methods that is based on the first approach 
is the method of auxiliary
sources (MAS) \cite{kaklamani2002, tsitsas2018}. 
The MAS (or close variants of MAS) is 
also known as the \emph{source model technique}
\cite{hochman2004, ludwig2004}, the \emph{method of fundamental solutions} \cite{Barnett2008, Cheng2020}, and the \emph{method of fictitious sources}
\cite{zolla1994, zolla1996}. 
In MAS, the equivalent surface current densities are defined on
fictitious surfaces, referred to as \emph{auxiliary surfaces}.
The pertinent SIEs are then discretized, with the equivalent current densities replaced by a number of \emph{auxiliary sources} (ASs).  
For the two-dimensional (2-D) geometries dealt with in this paper, the ASs are infinitely long electric (for transverse magnetic ($\rm{TM}$) polarization) or magnetic (for transverse electric ($\rm{TE}$) polarization) current filaments, which radiate cylindrical waves described by zero-order Hankel functions.  Despite its interesting features, the effectiveness and applicability of the MAS are restricted by open issues, including the placement of the ASs and the convergence of the fields. 
Regarding the convergence issues, it was shown in  \cite{fikioris2006, fikioris2007, fikioris2013, fikioris2015, Fikioris2018convergent} for perfectly conducting scatterers and in \cite{valagiannopoulos2012} for dielectric scatterers, that---under certain conditions on the placement of the auxiliary surface, and for sufficiently large numbers $N$ of ASs---the auxiliary currents diverge and oscillate while the generated MAS fields converge to the exact ones, provided that the oscillations are mild. However, for even larger values of $N$, the oscillations' amplitudes increase further and numerically corrupted fields are obtained.

A well-established method that utilizes the second approach is the \emph{Extended Boundary Condition Method} (EBCM) \cite{waterman1965, waterman1969} also known as \emph{Null Field Method} (NFM)  \cite{doicu20201, doicu20202w, wriedt2007,eremina2004, hellmers2008,egel2017}.
The EIEs are discretized, with the equivalent current densities on the physical boundaries replaced by $N$ discrete sources (similar to those of MAS). 
For problems involving perfectly conducting scatterers, it was shown that a main advantage of NFM over MAS is that in NFM no divergence and oscillations phenomena occur for the discrete currents \cite{fikioris2013, fikioris2011}.  
Particularly, for a circular cylinder, in \cite{fikioris2011} it was shown analytically that the EIE is solvable, irrespective of the position of the surface where the EBCs are applied, and that the solution of the discrete version of the EIE, when properly normalized, converges to the true solution as $N\rightarrow\infty$.
Hence, NFM showcases a very different behavior compared to MAS
despite the close resemblance of the integral equations involved in the two methods.

In this paper, we consider the problem of (internal or external) line-source excitation of a dielectric cylinder and combine the concepts of SEP and NFM by defining two equivalent problems of the original one and determining the unknown current densities on the physical boundary via a system of two integral equations stemming from two extended boundary conditions.
Our main goal is to study the convergence properties of the combined NFM-SEP approach and see how the above discussed conclusions on the solutions of the EIEs for perfectly conducting cylinders carry over to dielectric cylinders. As in \cite{fikioris2011}, we bring out stark contrasts to previous findings on the MAS, despite that the matrices of the NFM-SEP and the MAS are actually the same. 

The structure of the paper is as follows.
In Section \ref{sec:original}, we define the original scattering/excitation problems and apply the SEP to obtain its two equivalent problems.
Next, Section~\ref{sec:NFM-SEP-discrete} clarifies the implementation of the NFM-SEP to the considered problems. In Section~\ref{sec:circular_problem}, we give the exact solution for the circular problem and then consider and solve in Section~\ref{sec:continuous_EIE} the continuous version of the involved EIEs to derive analytically this exact solution. 
Then, in Section \ref{sec:Discrete_EIE}, we formulate the discrete version of the NFM-SEP and obtain the corresponding solutions. 
There, we prove that (as opposed to the MAS) no oscillations appear in the discrete NFM-SEP currents, which always converge to the corresponding continuous current densities.
Section~\ref{nfm-sep-vs-mas} contains a detailed comparative discussion between properties of  NFM-SEP and MAS. Furthermore, in Section \ref{sec:Numerical}, we show representative numerical results for circular and noncircular shapes illustrating some of the conclusions of Sections \ref{sec:continuous_EIE}-\ref{nfm-sep-vs-mas}. The paper is closed with conclusions in Section \ref{sec:conclusions}. 

An $e^{\mathrm i \omega t}$ dependence is assumed and suppressed in what follows, where $\omega$ is the angular frequency, assumed throughout to be real.


\section{NFM-SEP}
We describe NFM-SEP in the present section, beginning with the two equivalent problems.
\subsection{The Original Problem and its Equivalents}
\label{sec:original}

Figure \ref{fig:origina_equiv}a shows the cross section of an infinitely-long and parallel to the $z$-axis dielectric cylinder
(Region 2 or $R_2$), which is characterized by dielectric permittivity $\varepsilon_2$ and magnetic permeability $\mu_2$.
The cylinder is surrounded by a homogeneous dielectric medium, which is referred to as Region 1 ($R_1$) and characterized by respective parameters $\varepsilon_1$ and  $\mu_1$. The boundary between the two regions is denoted as $C$ and is of arbitrary shape. We consider two scattering problems:

\begin{figure}[htb!]
    \centering
    \subfigure[]{\includegraphics[width=0.52\textwidth]{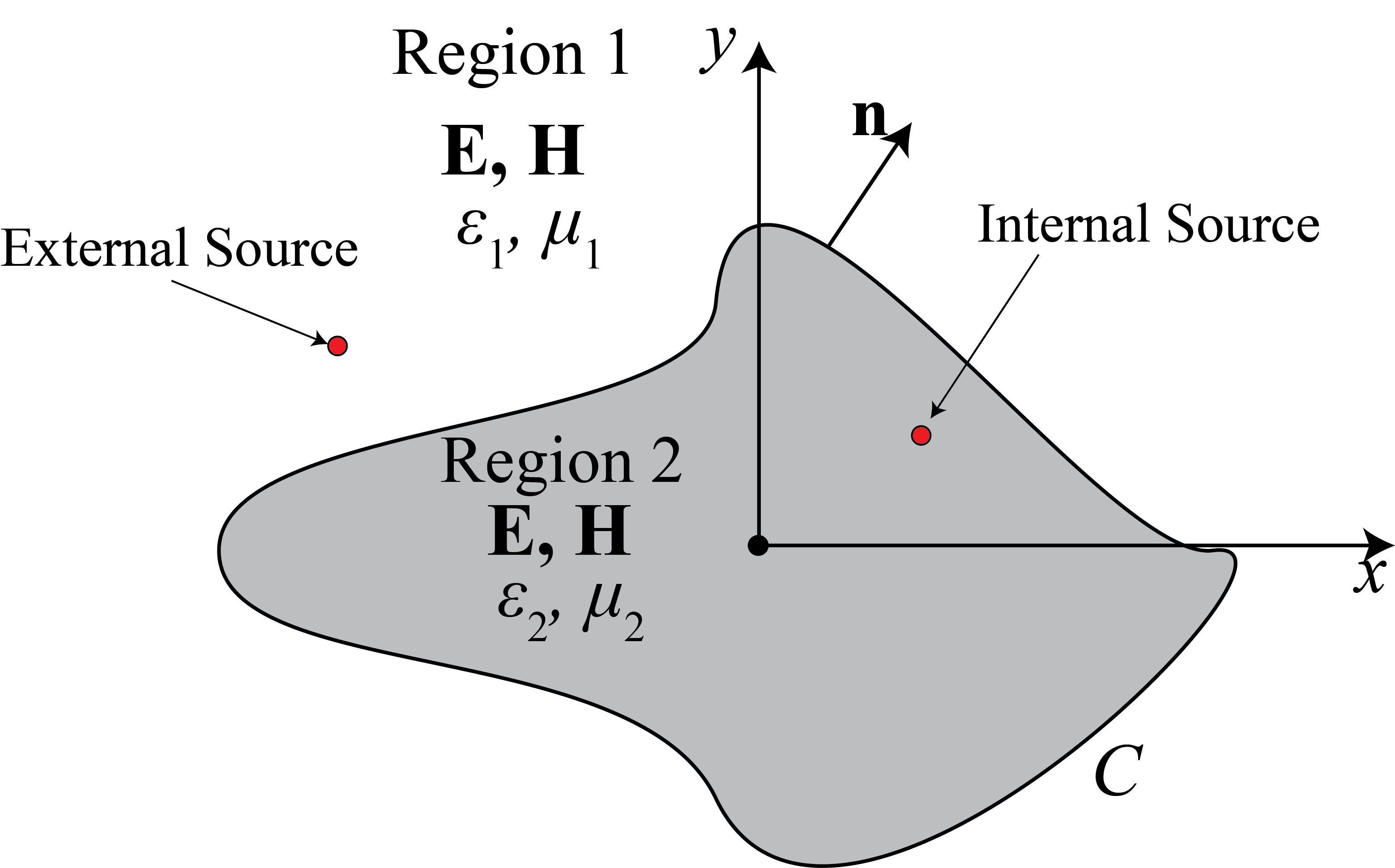}} 
    \subfigure[]{\includegraphics[width=0.40\textwidth]{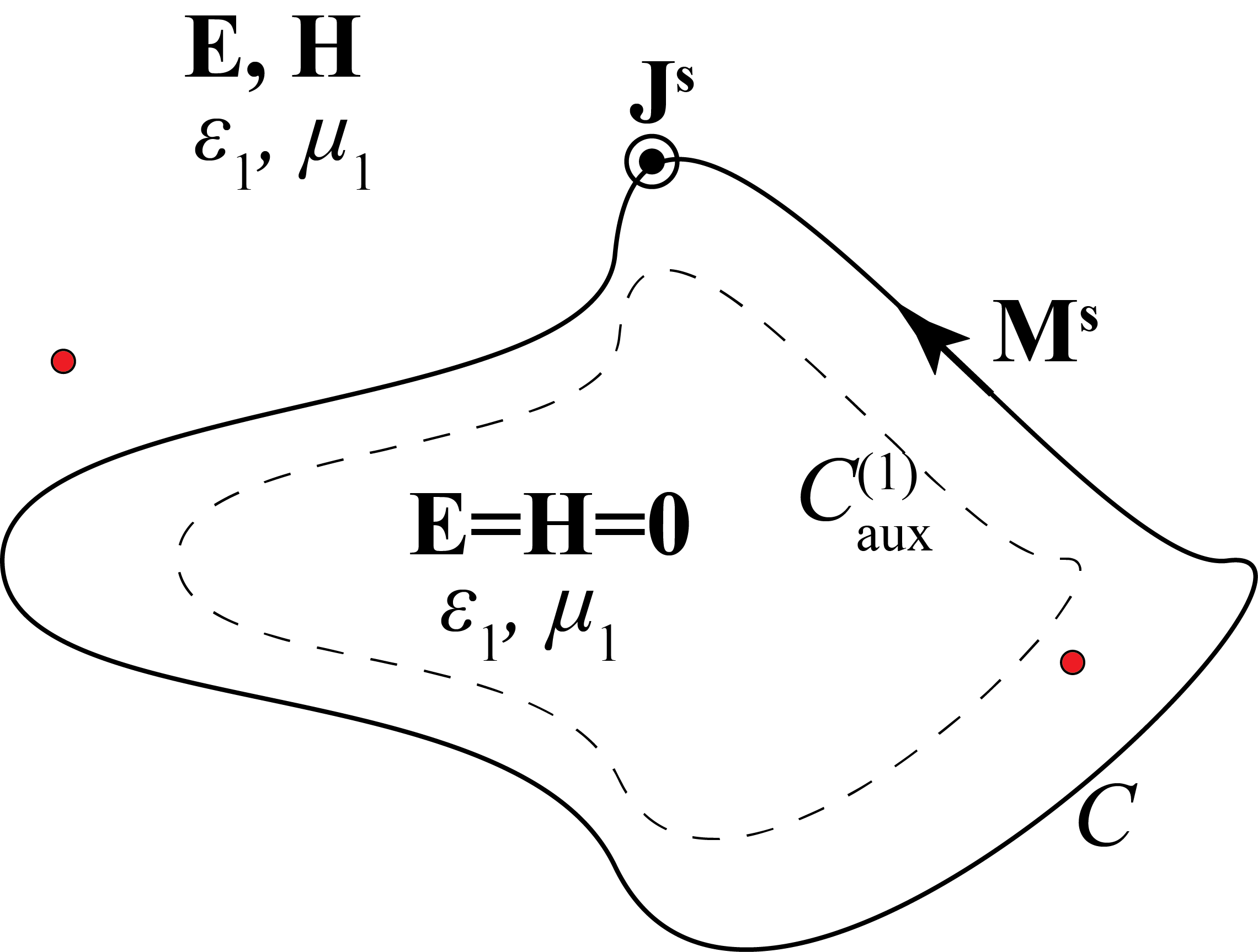}}
    \centering
    \subfigure[]{\includegraphics[width=0.60\textwidth]{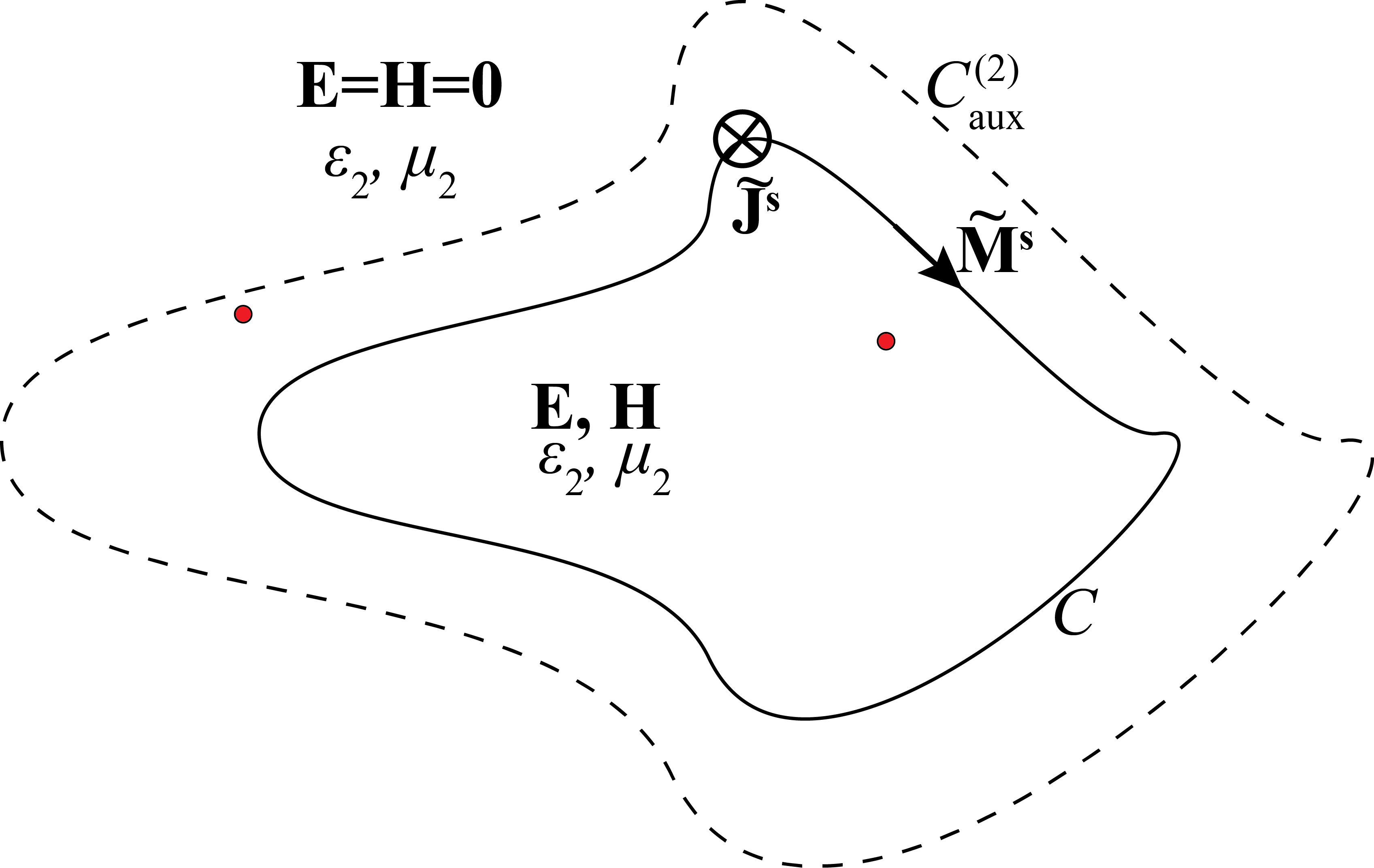}} 
    \caption{Scattering of a cylindrical wave by an infinitely-long dielectric cylinder. (a) Original problem; (b) first equivalent problem for the external field in $R_1$; and (c) second equivalent problem for the internal field in $R_2$.}
    \label{fig:origina_equiv}
\end{figure}

\begin{enumerate}
\item A $\rm{TM_z}$ or $\rm{TE_z}$ incident field propagates in $R_1$ and is scattered by the boundary $C$. The incident wave is either a plane wave (generated by a source far away from the cylinder), or a cylindrical wave generated by a $z$-directed electric ($\rm{TM_z}$) or magnetic ($\rm{TE_z}$) current filament of infinite length within $R_1$ and near $C$. This problem will be referred to as the \emph{external problem}.
\item An incident $\rm{TM_z}$ or $\rm{TE_z}$ field is generated by a $z$-directed electric or magnetic current filament of infinite length within $R_2$. This problem will be referred to as the \emph{internal problem}. 
\end{enumerate}

In both problems, the goal is to find the fields 
in $R_1$ and $R_2$. To do this, we apply the \emph{surface equivalence principle} (SEP) \cite{chew2022, Harrington_2001, jin2015}.
According to the SEP, two equivalent problems of the original one are defined, the first for the fields in $R_1$, and the second for the fields in $R_2$. 

In the first problem (depicted in Fig.~\ref{fig:origina_equiv}(b)), the fields in $R_1$ are the same as in the original problem, the material of $R_2$ is replaced by the material of $R_1$, and the fields in $R_2$ are set to zero. To satisfy the boundary conditions on $C$, two surface current densities, an electric $\mathbf J^s$ and a magnetic $\mathbf M^s$, are defined on $C$. These densities are given by 
\begin{align}
\label{eq:sources_1_J}
\mathbf J^{s}&=\hat{\mathbf n}\times \mathbf H,\\
\mathbf M^s&= \mathbf E \times \hat{\mathbf n},   
\label{eq:sources_1_M}
\end{align}
where $\hat{\mathbf n}$ is the outward normal unit vector on $C$, 
$\mathbf E$ and $\mathbf H$ is the total electric and magnetic field just outside $C$. In case of a $\rm{TM_z}$ ($\rm{TE_z}$) incident wave, $\mathbf J^s$ ($\mathbf M^s$) is $z$-directed, while $\mathbf M^s$ ($\mathbf J^s$) is tangential to $C$. 

In the second equivalent problem (depicted in Fig.~\ref{fig:origina_equiv}(c)), the fields in $R_2$ are the same as in the original problem, the material of $R_1$ is replaced by the material of $R_2$, and the fields in $R_1$ are set to zero. Now, we need to define the following current densities on $C$
\begin{align}
\label{eq:sources_2_J}
\Tilde{\mathbf J}^s&=-\hat{\mathbf n}\times \mathbf H,\\
\Tilde{\mathbf M}^s&=-\mathbf E \times \hat{\mathbf n},
\label{eq:sources_2_M}
\end{align}
i.e., the opposites of those defined in the first equivalent problem. 

For the external problem, the field in $R_1$ is the superposition of the incident field plus the field generated by the sources defined in (\ref{eq:sources_1_J}) and (\ref{eq:sources_1_M}), 
while the field in $R_2$ is that generated by the sources defined in (\ref{eq:sources_2_J}) and (\ref{eq:sources_2_M}). For the internal problem, the field in $R_2$ is the superposition of the incident field plus the field generated by the sources defined in (\ref{eq:sources_2_J}) and (\ref{eq:sources_2_M}) and the field in $R_1$ is that generated by the sources defined in (\ref{eq:sources_1_J}) and (\ref{eq:sources_1_M}). 

To implement the SEP, it is, thus, first necessary to determine 
$\mathbf J^s$ and $\mathbf M^s$. This can be accomplished by means of a system of two linear integral equations resulting from two extended boundary conditions.  In the context of NFM, the first such boundary condition arises from the first equivalent problem and requires the zeroing of the total (electric or magnetic) field on an auxiliary surface $C^{(1)}_{\mathrm{aux}}$ inside $R_2$. The second extended boundary condition stems from the second equivalent problem and requires the zeroing of the total field on an auxiliary surface $C^{(2)}_{\mathrm{aux}}$ inside $R_1$. 

The above-described ``continuous'' version of NFM-SEP leads to a system of linear integral equations. For brevity, we do not write this system down. Instead, we proceed to describe its discretized version, which is what one solves in actual practice. 

\subsection{Implementation of NFM-SEP}
\label{sec:NFM-SEP-discrete}

\begin{figure}[htb!]
    \centering
    \subfigure[]{\includegraphics[width=0.40\textwidth]{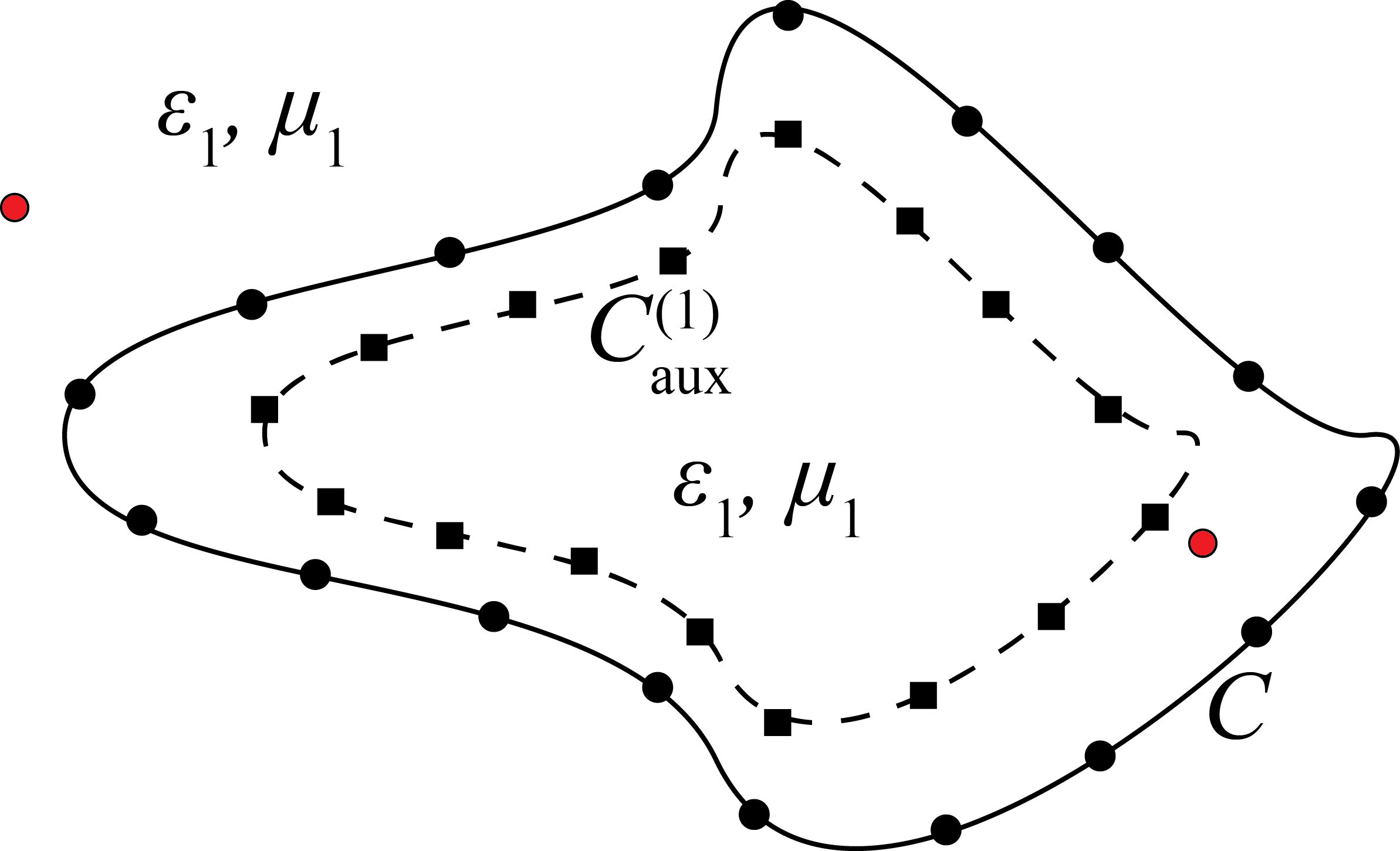}} 
    \subfigure[]{\includegraphics[width=0.51\textwidth]{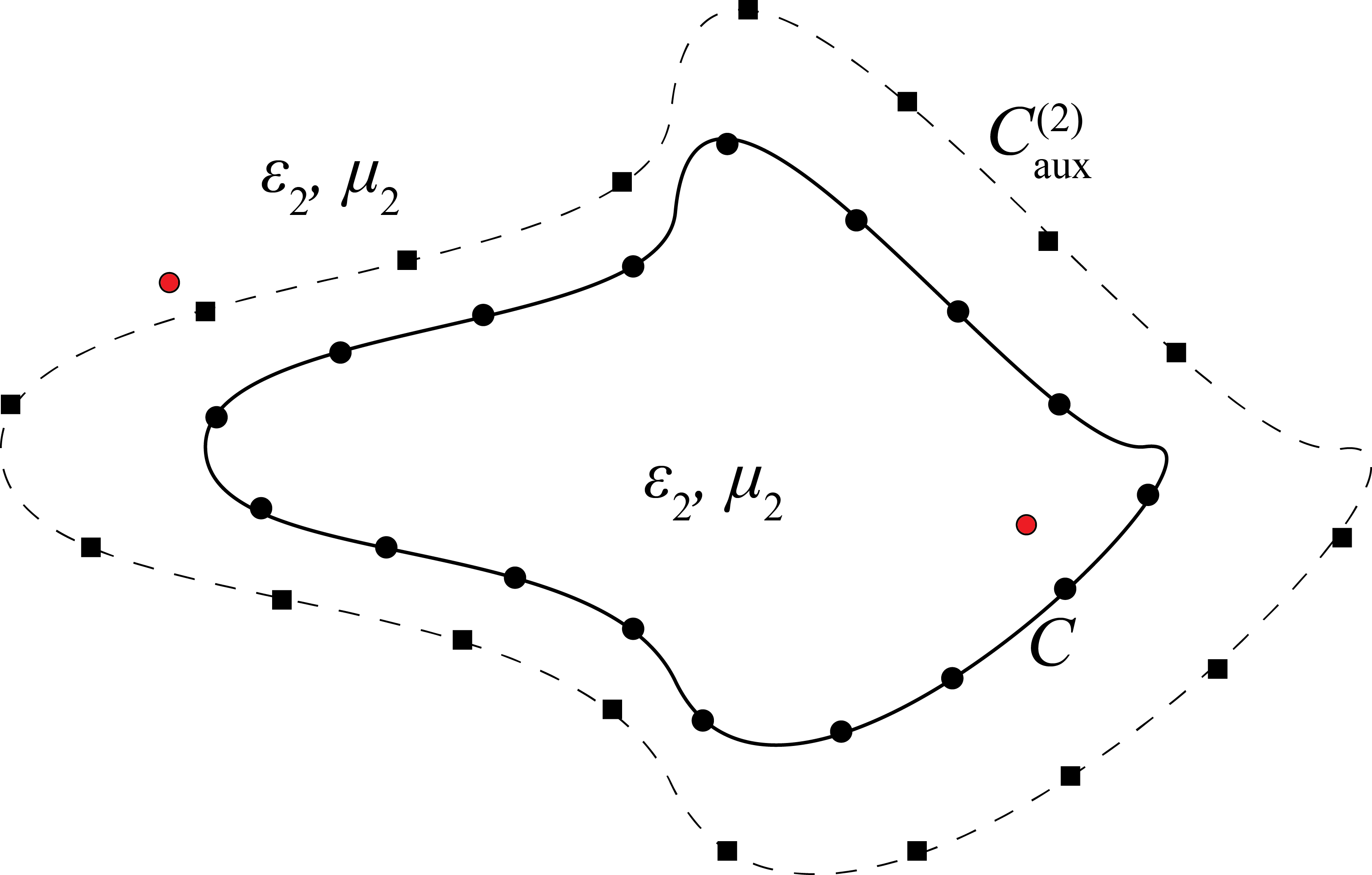}}
    \caption{Discretization of (a) the external and (b) the internal NFM equivalent problem for a dielectric cylinder of arbitrary cross section. Black dots represent the discrete electric and magnetic sources, and black squares represent the collocation points.}
    \label{fig:discrete_arb_NFM}
\end{figure}

The discretized versions of the NFM-SEP problems of Fig.~\ref{fig:origina_equiv} are depicted in Fig.~\ref{fig:discrete_arb_NFM}. Here, the continuous current densities are approximated by a number of discrete sources. Specifically, and in case of $\rm{TM_z}$ incidence, the $z$-directed electric current density $\mathbf J^s$
is replaced by $N$ infinitely-long and parallel to the $z$-axis electric current filaments with current complex amplitudes $I_l$. The tangential magnetic current density $\mathbf M^s$ is replaced by $N$ tangential on $C$ magnetic current filaments, each carrying a voltage $K_l$ ($l=0,1,\ldots,N-1$). 

The $z$-directed incident electric field is given by
\begin{equation}
 E^{\mathrm{inc}}_j(\mathbf r)=-\frac{k_jZ_j}{4} I H_0^{(2)}(k_j|\mathbf r-\mathbf r_{\mathrm{fil}}|),\,\,j=1,2,
 \label{eq:incident_electric}
\end{equation}
where $j$ denotes the region in which the filament is positioned, $I$ and $\mathbf r_{\mathrm{fil}}$ are the constant electric current amplitude and position vector of the filament,
$k_j=\omega\sqrt{\varepsilon_j\mu_j}$ and $Z_j=\sqrt{\mu_j/\varepsilon_j}$ are the wavenumber and impedance of region $j$, and  $H_0^{(2)}$ is the 0-order and second kind cylindrical Hankel function.

For external excitation, the total electric field in $R_1$ is given by 
\begin{align}
\nonumber
    E_z^{(1)}=&-\frac{k_1Z_1}{4}IH_0^{(2)}(k_1|\mathbf r-\mathbf r_{\mathrm{fil}}|)\\
    &-\frac{k_1Z_1}{4}\sum_{l=0}^{N-1}I_lH_0^{(2)}(k_1|\mathbf r-\mathbf r_l|)+ 
    \frac{k_1}{4\mathrm i}\sum_{l=0}^{N-1} K_l\frac{\partial H_0^{(2)}(k_1|\mathbf r-\mathbf r_l|)} {\partial n},
\label{eq:R1_field_discrete_arb}
\end{align}
while the total electric field in $R_2$ has the form
\begin{align}
E_z^{(2)}=-\frac{k_2Z_2}{4}\sum_{l=0}^{N-1}(-I_l)H_0^{(2)}(k_2|\mathbf r-\mathbf r_l|)+\frac{k_2}{4\mathrm i}\sum_{l=0}^{N-1} (-K_l)\frac{\partial H_0^{(2)}(k_2|\mathbf r-\mathbf r_l|)}{\partial n},
\label{eq:R2_field_discrete_arb}
\end{align}
where $\mathbf r_l$ is the position vector of the $l$-th collocated electric and magnetic discrete source.
Imposing the boundary conditions at $N$ collocation points on $C^{(1)}_{\mathrm {aux}}$ and $C^{(2)}_{\mathrm{aux}}$, we obtain the equations, for $p=0,1,\ldots,N-1$,
\begin{align}
\nonumber
&Z_1\sum_{l=0}^{N-1}I_lH_0^{(2)}(k_1|\mathbf r^{(1)}_{p,\mathrm{aux}}-\mathbf r_l|)+ \mathrm i\sum_{l=0}^{N-1} K_l\frac{\partial H_0^{(2)}(k_1|\mathbf r^{(1)}_{p,\mathrm{aux}}-\mathbf r_l|)} {\partial n}=\\
\label{eq:disc_linear_system_arb_1}
&-Z_1 I H_0^{(2)}(k_1|\mathbf r^{(1)}_{p,\mathrm{aux}}-\mathbf r_{\mathrm{fil}}|),\\
&Z_2\sum_{l=0}^{N-1}I_lH_0^{(2)}(k_2|\mathbf r^{(2)}_{p,\mathrm{aux}}-\mathbf r_l|)+\mathrm i\sum_{l=0}^{N-1}K_l\frac{\partial H_0^{(2)}(k_2|\mathbf r^{(2)}_{p,\mathrm{aux}}-\mathbf r_l|)} {\partial n}=0,
\label{eq:disc_linear_system_arb_2}
\end{align}
where $\mathbf r^{(j)}_{p,\mathrm{aux}}$ is the position vector of the $p$-th collocation point on $C^{(j)}_{\mathrm{aux}}$. Equations (\ref{eq:disc_linear_system_arb_1}) and (\ref{eq:disc_linear_system_arb_2}) constitute an $2N \times 2N$ linear system with unknowns the electric $I_l$ and magnetic $K_l$ current complex amplitudes. Once solved, the fields in each region are determined from (\ref{eq:R1_field_discrete_arb}) and (\ref{eq:R2_field_discrete_arb}).

\section{Circular Problem: Exact Solution and its Analytic Continuation}
\label{sec:circular_problem}

\begin{figure}[htb!]
    \centering
    \subfigure[]{\includegraphics[width=0.47\textwidth]{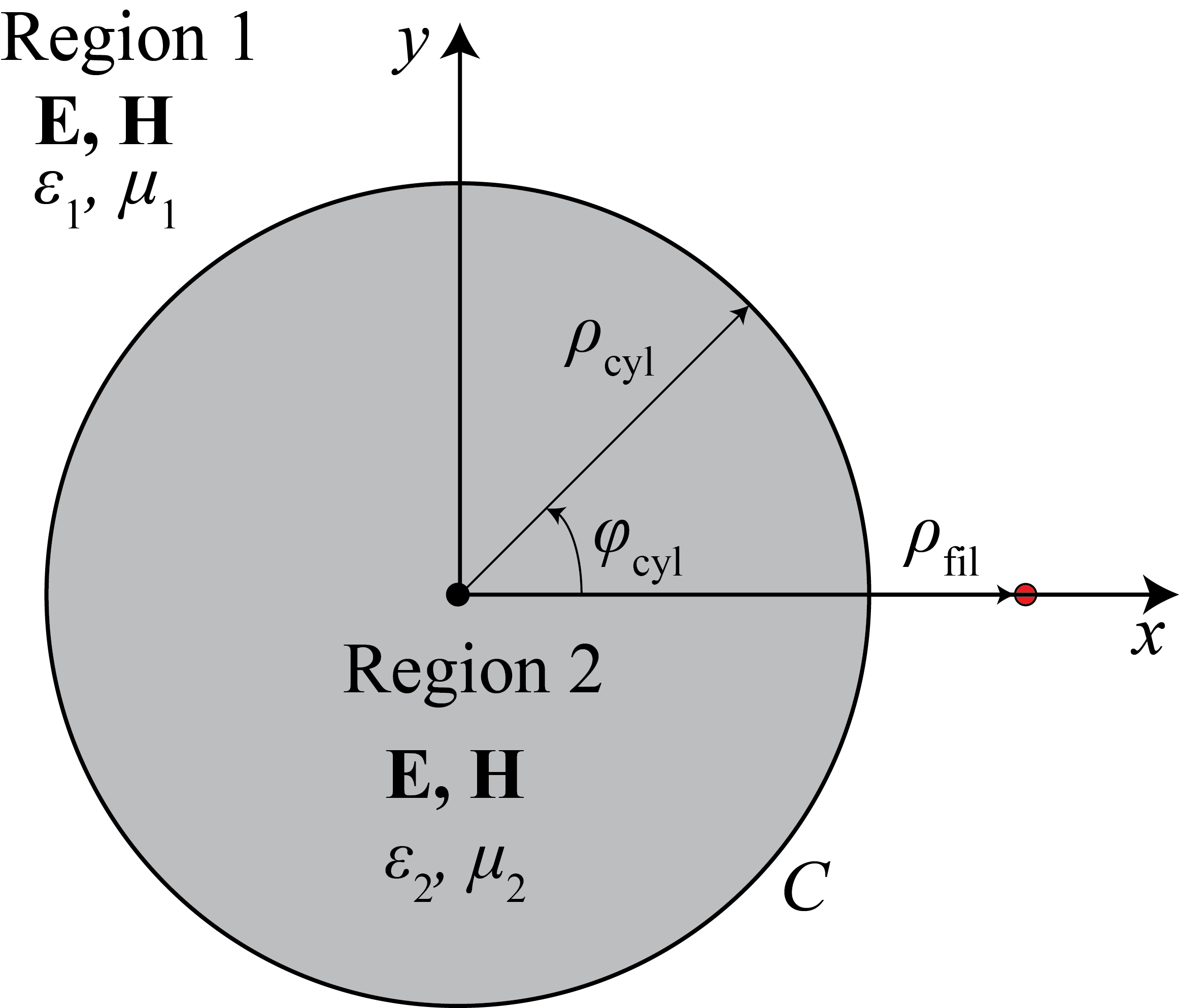}} 
    \subfigure[]{\includegraphics[width=0.47\textwidth]{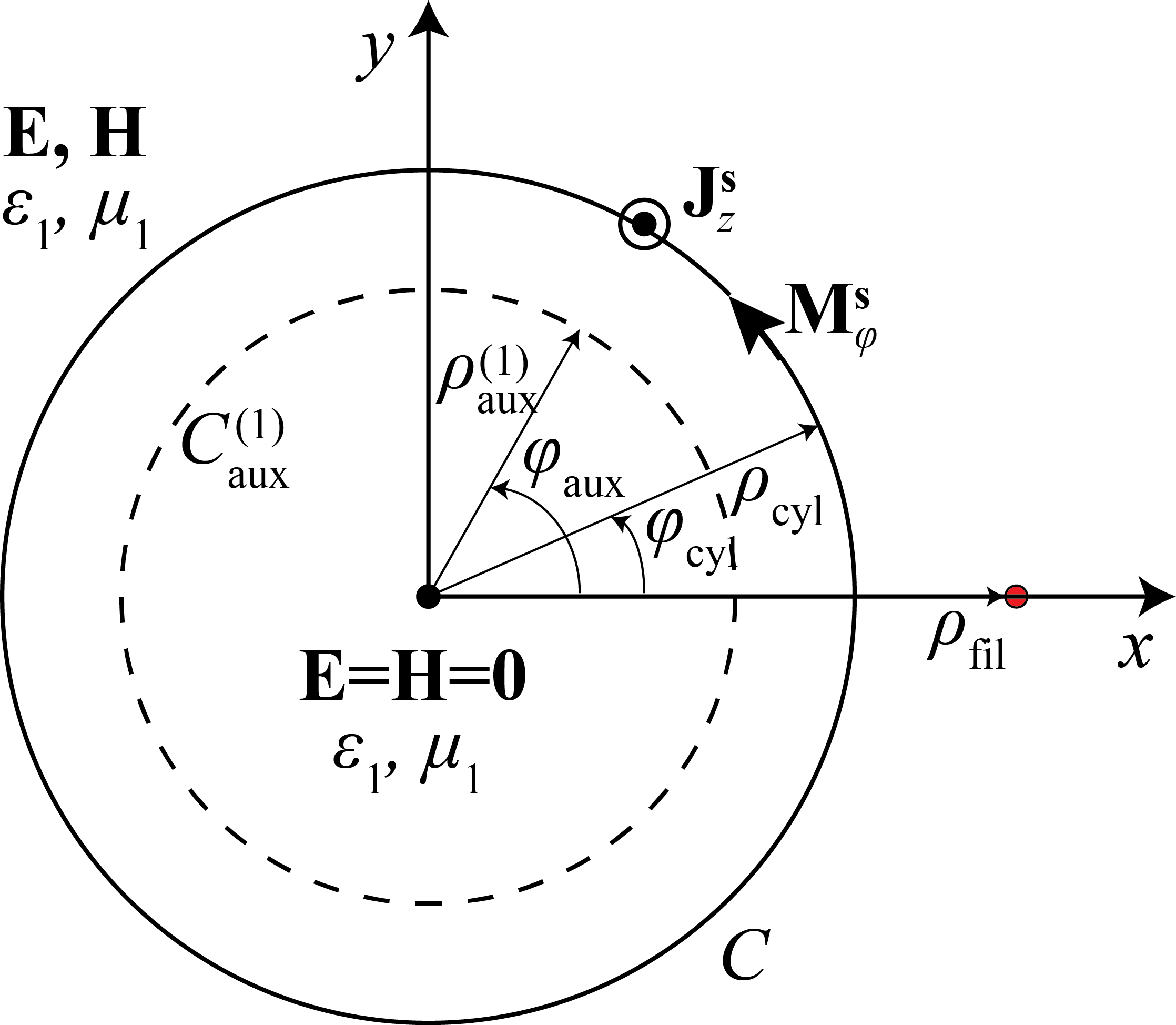}}
    \centering
    \subfigure[]{\includegraphics[width=0.52\textwidth]{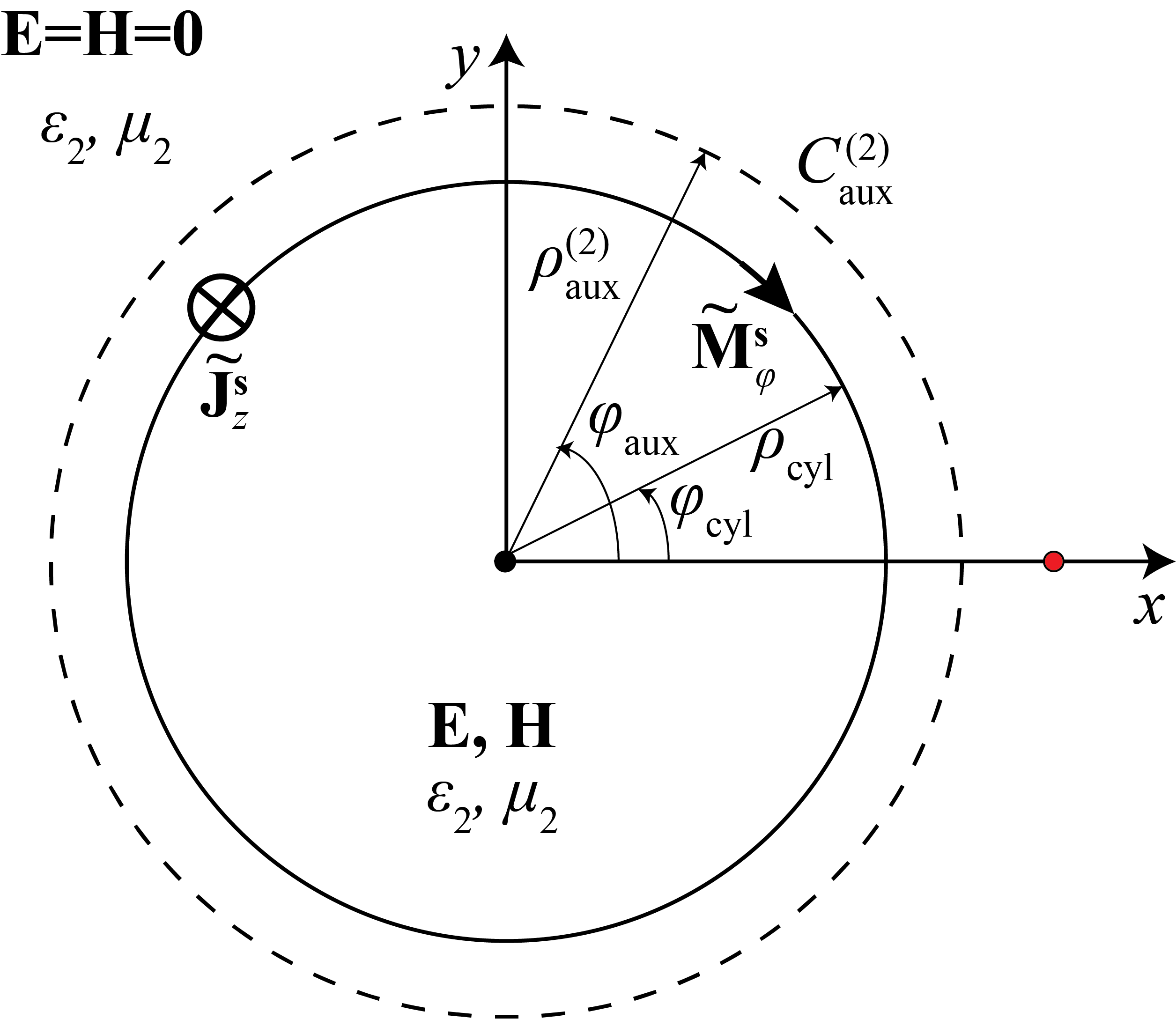}} 
    \caption{Scattering of a cylindrical wave by an infinitely long circular dielectric cylinder. (a) Original problem; (b) first equivalent problem (for external field in $R_1$); and (c) second equivalent problem (for internal field in $R_2$).}
    \label{fig:origina_equiv_circ}
\end{figure}

The simplest geometry is that of a circular boundary $C$, as illustrated in Fig.~\ref{fig:origina_equiv_circ}(a): A circular cylinder of radius $\rho_{\mathrm{cyl}}$ is excited by an external or an internal (electric or magnetic) source lying at $(\rho_{\mathrm{fil}}, 0)$ in  polar coordinates. For $\rm{TM_z}$ incidence, the sole $z$-component of the incident electric field is
\begin{equation}
 E_j^{\mathrm{inc}}(\rho_{\mathrm{obs}},\phi_{\mathrm{obs}})=-\frac{k_jZ_j}{4}IH_0^{(2)}(k_jD_{\mathrm{fil,obs}}),\,\,j=1,2,
 \label{eq:incident_electric}
\end{equation}
with $j$ the region in which the filament of amplitude $I$ is located, and $D_{\mathrm{fil,obs}}$ the distance between the filament position and the observation point, i.e., 
\begin{equation}
D_{\mathrm{fil,obs}}=\sqrt{\rho^2_{\mathrm{obs}}+\rho^2_{\mathrm{fil}}-2\rho_{\mathrm{obs}}\rho_{\mathrm{fil}}\cos{\phi_{\mathrm{obs}}}}.
\end{equation}

By the duality principle \cite{Harrington_2001, jin2015}, the incident $\rm TE_z$ magnetic field is given by
\begin{equation}
 H_j^{\mathrm{inc}}(\rho_{\mathrm{obs}},\phi_{\mathrm{obs}})=-\frac{k_j}{4Z_j}KH_0^{(2)}(k_jD_{\mathrm{fil,obs}}),\,\,j=1,2,
 \label{eq:incident_magnetic}
\end{equation}
where $K$ is the magnetic current amplitude of the filament.

The analysis that follows concerns $\rm{TM_z}$ polarization. The corresponding results for $\rm{TE_z}$ polarization are omitted for brevity, but can be easily obtained via the duality principle.

The problem of Fig.~\ref{fig:origina_equiv_circ}a is solved analytically by Fourier series \cite{valagiannopoulos2012}. For external excitation, the total electric fields in $R_1$ and $R_2$ are found to be
\begin{align}
\nonumber
&E^{(1)}_z(\rho_{\mathrm{obs}},\phi_{\mathrm{obs}})=E_1^{\mathrm{inc}}(\rho_{\mathrm{obs}},\phi_{\mathrm{obs}})+\frac{k_1Z_1I}{4}\sum_{n=-\infty}^{+\infty}H_n^{(2)}(k_1\rho_{\mathrm{obs}})H_n^{(2)}(k_1\rho_{\mathrm{fil}})\times \\&\frac{Z_1J'_n(k_2\rho_{\mathrm{cyl}})J_n(k_1\rho_{\mathrm{cyl}})-Z_2J_n(k_2\rho_{\mathrm{cyl}})J'_n(k_1\rho_{\mathrm{cyl}})}{Z_1H_n^{(2)}(k_1\rho_{\mathrm{cyl}})J'_n(k_2\rho_{\mathrm{cyl}})-Z_2J_n(k_2\rho_{\mathrm{cyl}})H'^{(2)}_n(k_1\rho_{\mathrm{cyl}})}e^{\mathrm i n \phi_{\mathrm{obs}}},
\label{eq:exact_ext_TM_1}
\\
\nonumber
&E^{(2)}_z(\rho_{\mathrm{obs}},\phi_{\mathrm{obs}})=-\frac{I}{2\pi\rho_{\mathrm{cyl}}}\sum_{n=-\infty}^{+\infty}J_n(k_2\rho_{\mathrm{obs}})H_n^{(2)}(k_1\rho_{\mathrm{fil}})\times \\
&\frac{\mathrm iZ_1Z_2}{Z_1H_n^{(2)}(k_1\rho_{\mathrm{cyl}})J'_n(k_2\rho_{\mathrm{cyl}})-Z_2J_n(k_2\rho_{\mathrm{cyl}})H'^{(2)}_n(k_1\rho_{\mathrm{cyl}})}e^{\mathrm i n \phi_{\mathrm{obs}}}.
\label{eq:exact_ext_TM_2}
\end{align}

For internal excitation, the respective total fields are
\begin{align}
\nonumber
&E^{(1)}_z(\rho_{\mathrm{obs}},\phi_{\mathrm{obs}})=-\frac{I}{2\pi\rho_{\mathrm{cyl}}}\sum_{n=-\infty}^{+\infty}H^{(2)}_n(k_1\rho_{\mathrm{obs}})J_n(k_2\rho_{\mathrm{fil}})\times \\
&\frac{\mathrm iZ_1Z_2}{Z_1H^{(2)}_n(k_1\rho_{\mathrm{cyl}})J'_n(k_2\rho_{\mathrm{cyl}})-Z_2J_n(k_2\rho_{\mathrm{cyl}})H'^{(2)}_n(k_1\rho_{\mathrm{cyl}})}e^{\mathrm i n \phi_{\mathrm{obs}}},
\label{eq:exact_int_TM_1}
\\
\nonumber
&E^{(2)}_z(\rho_{\mathrm{obs}},\phi_{\mathrm{obs}})=E_2^{\mathrm{inc}}(\rho_{\mathrm{obs}},\phi_{\mathrm{obs}})+\frac{k_2Z_2I}{4}\sum_{n=-\infty}^{+\infty}J_n(k_2\rho_{\mathrm{obs}})J_n(k_2\rho_{\mathrm{fil}})\times \\
&\frac{Z_1H^{(2)}_n(k_1\rho_{\mathrm{cyl}})H'^{(2)}_n(k_2\rho_{\mathrm{cyl}})-Z_2H'^{(2)}_n(k_1\rho_{\mathrm{cyl}})H^{(2)}_n(k_2\rho_{\mathrm{cyl}})}{Z_1H^{(2)}_n(k_1\rho_{\mathrm{cyl}})J'_n(k_2\rho_{\mathrm{cyl}})-Z_2J_n(k_2\rho_{\mathrm{cyl}})H'^{(2)}_n(k_1\rho_{\mathrm{cyl}})}e^{\mathrm i n \phi_{\mathrm{obs}}}.
\label{eq:exact_int_TM_2}
\end{align}

Note that in each term of all series (\ref{eq:exact_ext_TM_1})-(\ref{eq:exact_int_TM_2}) the denominator is the same. This denominator does not vanish for any real frequency value  \cite{marcuse1982, balanis2012}, and, hence,
we cannot detect the cutoff frequency of any propagating mode in the dielectric cylinder. This is in contrast to the case of PEC waveguides excited by an internal or external source \cite{kou_tsi_fik_21}, where the denominator of the expression of the total field within the PEC cavity vanishes for real frequency values corresponding to the cutoff frequencies. This is because in dielectric waveguides, the concept of cutoff is related to the inability to confine a propagating wave near the core \cite{marcuse1982}, whereas in PEC waveguides, to the inability of a propagating wave to exist. In other words, to detect cutoff frequencies in a dielectric cylinder, we would need to excite it with sources that induce some propagating wave along the axis of the cylinder. This cannot be achieved with an electric or magnetic filament with a constant amplitude as a source, but rather with a proper combination of both \cite{kou_tsi_fik_22}, or by exciting the dielectric cylinder with an oblique plane wave \cite{wait1955}.

Needless to say, the four series converge within their region of validity, i.e., the two series associated with $E_z^{(1)}$ ($E_z^{(2)}$) converge for $\rho_\mathrm{obs}>\rho_\mathrm{cyl}$~($\rho_\mathrm{obs}<\rho_\mathrm{cyl}$). But, we can show that all series converge within larger regions; and that each function, thus, defined is in fact the analytic continuation of the respective field into the extended region. The four extended regions have circular boundaries and, as seen in Table~\ref{table:convergence}, are described by the relation between $\rho_\mathrm{obs}$ on the one hand, and either $\rho_\mathrm{fil}$ or the ``critical  distance'' $\rho_{\mathrm{cri}}$ on the other. This distance is defined by
\begin{equation}
\rho_{\mathrm{cri}}=\frac{\rho_{\mathrm{cyl}}^2}{\rho_{\mathrm{fil}}}.
\label{eq:critical}
\end{equation}
\begin{table}[]
\caption{Each cell shows the extended region of convergence of a series in (\ref{eq:exact_ext_TM_1}), (\ref{eq:exact_ext_TM_2}) or (\ref{eq:exact_int_TM_1}), (\ref{eq:exact_int_TM_2}).  
For example, the series in (\ref{eq:exact_int_TM_1}) pertains to internal excitation (so that $\rfil<\rcyl$) and gives the field in $R_1$ (i.e., for observation points $\robs>\rcyl$); but this series converges not only in $\robs>\rcyl$, but in the extended region $\robs>\rfil$. All series diverge in the complementary regions; for example, the series in (\ref{eq:exact_int_TM_1}) diverges for $\robs<\rfil$.}
\label{table:convergence}
\centering
\begin{tabular}{c|c|c|}
\cline{2-3}
{}                                  & $R_2,\,\robs<\rcyl$  & $R_1,\,\robs>\rcyl$  \\ \hline
\multicolumn{1}{|c|} {external excitation,$\,\rcyl<\rfil$} & (\ref{eq:exact_ext_TM_2}),$\,\robs<\rfil$ & (\ref{eq:exact_ext_TM_1}),$\,\robs>\rcri$ \\ \hline
\multicolumn{1}{|c|} {internal excitation,$\,\rfil<\rcyl$} & (\ref{eq:exact_int_TM_2}),$\,\robs<\rcri$ & (\ref{eq:exact_int_TM_1}),$\,\robs>\rfil$ \\ \hline
\end{tabular}
\end{table}
To derive the results in Table \ref{table:convergence}, use the large-order asymptotic approximations (for the Bessel/Hankel functions and their derivatives) provided in  \ref{sec:appendix}. 
For example, when $|n|$ is large, the $n$-th terms of the series in (\ref{eq:exact_ext_TM_1}) (scattered field) and (\ref{eq:exact_ext_TM_2}) are seen to behave like
\begin{align}
\label{eq:convergence_exact_1}
&\frac{2}{\pi|n|}\frac{\rho_{\mathrm{cri}}^{|n|}}{\rho_{\mathrm{obs}}^{|n|}}e^{\mathrm i n \phi_{\mathrm{obs}}},\\
&\frac{\mathrm i Z_1Z_2}{\frac{Z_1}{k_2}-\frac{Z_2}{k_1}}\frac{\rho_{\mathrm{obs}}^{|n|}}{\rho_{\mathrm{fil}}^{|n|}}e^{\mathrm i n \phi_{\mathrm{obs}}},
\label{eq:convergence_exact_2}
\end{align}
respectively, where we used the definition (\ref{eq:critical}). 
By (\ref{eq:convergence_exact_1}), the series pertaining to $R_1$ converges exponentially for $\rho_{\mathrm{obs}}>\rho_{\mathrm{cri}}$ and diverges exponentially for $\rho_{\mathrm{obs}}<\rho_{\mathrm{cri}}$, verifying the pertinent entry of Table~1.
Since the derivative of (\ref{eq:exact_ext_TM_1}) is continuous across $\robs=\rho_{\mathrm{cyl}}$, this extended field is in fact the analytic continuation of the scattered field (series in (\ref{eq:exact_ext_TM_1}))  from the original region $R_1$ into \textit{part} of the complementary region $R_2$. 

Recall that, in MAS, an ``auxiliary surface'' is a surface on which  MAS (auxiliary) currents are located. 
In the case of MAS, it is found in \cite{valagiannopoulos2012} that placing an auxiliary surface beyond the respective extended region, as defined in Table~\ref{table:convergence}, results in divergence and unphysical oscillations in the corresponding auxiliary currents. 
The conclusions regarding the divergence of the auxiliary currents are summarized in Table~\ref{table:divergence-MAS-currents}.
\begin{table}[]
\caption{Each cell shows the region of divergence of the MAS auxiliary currents. For example, in case of internal excitation (so that $\rfil<\rcyl$), the MAS currents on the internal (external) auxiliary surface $C^{(1)}_{\mathrm{aux}}$ $(C^{(2)}_{\mathrm{aux}})$ diverge when $\,\rho^{(1)}_{\mathrm{aux}}<\rfil$ ($\,\rho^{(2)}_{\mathrm{aux}}>\rcri$).}
\label{table:divergence-MAS-currents}
\centering
\begin{tabular}{c|c|c|}
\cline{2-3}
{}                                  & $C^{(1)}_{\mathrm{aux}}$  & $C^{(2)}_{\mathrm{aux}}$  \\ \hline
\multicolumn{1}{|c|} {external excitation,$\,\rcyl<\rfil$} & $\,\rho^{(1)}_{\mathrm{aux}}<\rcri$ & $\,\rho^{(2)}_{\mathrm{aux}}>\rfil$ \\ \hline
\multicolumn{1}{|c|} {internal excitation,$\,\rfil<\rcyl$} & $\,\rho^{(1)}_{\mathrm{aux}}<\rfil$ & $\,\rho^{(2)}_{\mathrm{aux}}>\rcri$ \\ \hline
\end{tabular}
\end{table}

In the next section, we show that the above mentioned findings concerning the divergence and oscillations in the MAS auxiliary currents have no counterpart in the case of the EIE.

\section{Continuous EIE Solution and its Convergence}
\label{sec:continuous_EIE}

The key idea is to analytically extract the solution of the circular problem by applying the NFM. The two equivalent circular problems are depicted in Figs.~\ref{fig:origina_equiv_circ}(b) and (c) respectively. In case of external $\rm{TM_z}$ incidence, according to the first equivalent problem of Fig.~\ref{fig:origina_equiv_circ}b, the total electric field in $R_1$ is the following superposition of the incident field and the field radiated by the current densities $J_z^s$ and $M_\phi^s$
\begin{align}
\nonumber
&E_z^{(1)}(\rho_{\mathrm{obs}}, \phi_{\mathrm{obs}})=
E_1^{\mathrm{inc}}(\rho_{\mathrm{obs}}, \phi_{\mathrm{obs}})-\frac{k_1Z_1\rho_{\mathrm {cyl}}}{4}\int_{-\pi}^{\pi}J^s_z(\phi_{\mathrm{cyl}})H_0^{(2)}(k_1D_{\mathrm{obs,cyl}})\mathrm{d}\phi_{\mathrm{cyl}}\\
&-\frac{k_1\rho_{\mathrm{cyl}}}{4\mathrm i}\int_{-\pi}^{\pi}M^s_\phi(\phi_{\mathrm{cyl}})\frac{\partial H_0^{(2)}(k_1D_{\mathrm{obs,cyl}}) }{\partial \rcyl}\mathrm{d}\phi_{\mathrm{cyl}},
\label{eq:R1_field}
\end{align}
with
\begin{equation}
D_{\mathrm{obs,cyl}}=\sqrt{\rho^2_{\mathrm{cyl}}+\rho^2_{\mathrm{obs}}-2\rho_{\mathrm{cyl}}\rho_{\mathrm{obs}}\cos\mathrm ({\phi_{\mathrm{cyl}}}-\phi_{\mathrm{obs}})}
\label{eq:cosine_obs_cyl}
\end{equation}
the distance between a point on the circular boundary $C$ and the observation point. 
According to the second equivalent problem (Fig.~\ref{fig:origina_equiv_circ}c), the field in $R_2$ is the following superposition of the fields radiated by $-J_z^s$ and $-M_\phi^s$ in an unbounded medium filled with the cylinder's material
\begin{align}
\nonumber
 E_z^{(2)}(\rho_{\mathrm{obs}}, \phi_{\mathrm{obs}})= 
 &-\frac{k_2Z_2\rho_{\mathrm {cyl}}}{4}\int_{-\pi}^{\pi}(-J^s_z(\phi_{\mathrm{cyl}}))H_0^{(2)}(k_2D_{\mathrm{obs,cyl}})\mathrm{d}\phi_{\mathrm{cyl}}\\&-\frac{k_2\rho_{\mathrm{cyl}}}{4\mathrm i}\int_{-\pi}^{\pi}(-M^s_\phi(\phi_{\mathrm{cyl}}))\frac{\partial H_0^{(2)}(k_2D_{\mathrm{obs,cyl}}) }{\partial \rcyl}\mathrm{d}\phi_{\mathrm{cyl}}.
 \label{eq:R2_field}
\end{align}

The unknown current densities $J_z^s$ and $M_\phi^s$ in Eqs.~(\ref{eq:R1_field}) and (\ref{eq:R2_field}) are determined by enforcing the two boundary conditions discussed in Section \ref{sec:original}. The first stems from the first equivalent problem and requires the zeroing of the total field on the auxiliary surface $C^{(1)}_{\mathrm{aux}}$, with $\rho^{(1)}_{\mathrm{aux}}<\rho_{\mathrm{cyl}}$, yielding
\begin{align}
\nonumber
&-\frac{k_1Z_1\rho_{\mathrm {cyl}}}{4}\int_{-\pi}^{\pi}J^s_z(\phi_{\mathrm{cyl}})H_0^{(2)}(k_1D^{(1)}_{\mathrm{cyl,aux}})\mathrm{d}\phi_{\mathrm{cyl}}\\
&+
\nonumber
\frac{k_1\rho_{\mathrm{cyl}}}{4\mathrm i}\int_{-\pi}^{\pi}M^s_\phi(\phi_{\mathrm{cyl}})\frac{\rho_{\mathrm{cyl}}-\rho^{(1)}_{\mathrm{aux}}\cos(\phi_{\mathrm{cyl}}-\phi_{\mathrm{aux}})}{D^{(1)}_{\mathrm{cyl,aux}}}H_1^{(2)}(k_1D^{(1)}_{\mathrm{cyl,aux}})\mathrm{d}\phi_{\mathrm {cyl}}=\\
&
\frac{k_1Z_1}{4}IH_0^{(2)}(k_1 D^{(1)}_{\mathrm{fil,aux}}),
\label{eq:cont_eq_1}
\end{align}
where $D^{(1)}_{\mathrm{fil,aux}}$ is the distance between the filament $I$ and a point on $C^{(1)}_{\mathrm{aux}}$ with coordinates $(\rho^{(1)}_{\mathrm{aux}},\phi_{\mathrm{aux}})$, i.e.,
\begin{equation}
D^{(1)}_{\mathrm{fil, aux}}=\sqrt{\rho^2_{\mathrm{fil}}+(\rho^{(1)}_{\mathrm{aux}})^2-2\rho_{\mathrm{fil}}\rho^{(1)}_{\mathrm{aux}}\cos\phi^{(1)}_{\mathrm{aux}}},
\end{equation}
and $D^{(1)}_{\mathrm{cyl,aux}}$ is the distance from a point on $C$ to a point on $C^{(1)}_{\mathrm{aux}}$
\begin{equation}
D^{(1)}_{\mathrm{cyl, aux}}=\sqrt{\rho^2_{\mathrm{cyl}}+(\rho^{(1)}_{\mathrm{aux}})^2-2\rho_{\mathrm{cyl}}\rho^{(1)}_{\mathrm{aux}}\cos(\mathrm {\phi_{\mathrm{cyl}}}-\phi^{(1)}_{\mathrm{aux}})}.
\end{equation}

From the second equivalent problem, and for $\rho^{(2)}_{\mathrm{aux}}>\rho_{\mathrm{cyl}}$, we obtain
\begin{align}
\nonumber
&-\frac{k_2Z_2\rho_{\mathrm {cyl}}}{4}\int_{-\pi}^{\pi}J^s_z(\phi_{\mathrm{cyl}})H_0^{(2)}(k_2D^{(2)}_{\mathrm{cyl,aux}})\mathrm{d}\phi_{\mathrm{cyl}}\\&+\frac{k_2\rho_{\mathrm{cyl}}}{4\mathrm i}\int_{-\pi}^{\pi}M^s_\phi(\phi_{\mathrm{cyl}})\frac{\rho_{\mathrm{cyl}}-\rho^{(2)}_{\mathrm{aux}}\cos(\phi_{\mathrm{cyl}}-\phi^{(2)}_{\mathrm{aux}})}{D^{(2)}_{\mathrm{cyl,aux}}}H_1^{(2)}(k_2D^{(2)}_{\mathrm{cyl,aux}})\mathrm{d}\phi_{\mathrm {cyl}}=0,
\label{eq:cont_eq_2}
\end{align}
where 
\begin{equation}
D^{(2)}_{\mathrm{cyl, aux}}=\sqrt{\rho^2_{\mathrm{cyl}}+(\rho^{(2)}_{\mathrm{aux}})^2-2\rho_{\mathrm{cyl}}\rho^{(2)}_{\mathrm{aux}}\cos(\mathrm {\phi_{\mathrm{cyl}}}-\phi^{(2)}_{\mathrm{aux}})}.
\end{equation}

Together, (\ref{eq:cont_eq_1}) and (\ref{eq:cont_eq_2}) consist a $2\times 2$ system of Fredholm linear integral equations, in which the involved kernels are even and $2\pi$-periodic. Such systems
can be solved explicitly by using Fourier series (as described, e.g.,  in Appendix A of \cite{fikioris2006}). Let
 \begin{align}
J^s_z(\phi_{\mathrm{cyl}})&=\sum_{n=-\infty}^{+\infty}I_n^s e^{\mathrm i n \phi_{\mathrm {cyl}}},
\label{eq:Fourier_series_1}
\\
M^s_\phi(\phi_{\mathrm{cyl}})&=\sum_{n=-\infty}^{+\infty}M_n^s e^{\mathrm i n \phi_{\mathrm {cyl}}},
 \label{eq:Fourier_series_2}
 \end{align}
where $I_n^s$ and  $M_n^s$ are the respective Fourier coefficients. 
Then, we multiply Eqs.~(\ref{eq:cont_eq_1}) and (\ref{eq:cont_eq_2}) with $e^{\mathrm i n \phi^{(j)}_{\mathrm{aux}}}$, and integrate with respect to $\phi^{(j)}_{\mathrm{aux}}$, with $j=1,2$. Next, we interchange the orders of integration and invoke the addition theorems (\ref{eq:addition_theorem})-(\ref{eq:addition_theorem_der_2}) to introduce the Fourier coefficients
\begin{align}
    &H_n^{(2)}(k_1\rho_{\mathrm{cyl}})J_n(k_1\rho^{(1)}_{\mathrm{aux}})e^{-\mathrm i n \phi_{\mathrm{cyl}}}=\frac{1}{2\pi}\int_{-\pi}^{\pi}H_0^{(2)}(k_1D^{(1)}_{\mathrm{cyl,aux}})e^{-\mathrm i n \phi^{(1)}_{\mathrm{aux}}}\mathrm{d}\phi^{(1)}_{\mathrm{aux}}\\
 &H_n^{(2)}(k_2\rho^{(2)}_{\mathrm{aux}})J_n(k_2\rho_{\mathrm{cyl}})e^{-\mathrm i n \phi_{\mathrm{cyl}}}=\frac{1}{2\pi}\int_{-\pi}^{\pi}H_0^{(2)}(k_2D^{(2)}_{\mathrm{cyl,aux}})e^{-\mathrm i n \phi^{(2)}_{\mathrm{aux}}}\mathrm{d}\phi^{(2)}_{\mathrm{aux}}\\
 &H_n^{'(2)}(k_1\rho_{\mathrm{cyl}})J_n(k_1\rho^{(1)}_{\mathrm{aux}})e^{-\mathrm i n \phi_{\mathrm{cyl}}}=\\ 
 &-\frac{1}{2\pi}\int_{-\pi}^{\pi}\frac{\rho_{\mathrm{cyl}}-\rho^{(1)}_{\mathrm{aux}}\cos(\phi_{\mathrm{cyl}}-\phi^{(1)}_{\mathrm{aux}})}{D^{(1)}_{\mathrm{cyl,aux}}}H_1^{(2)}(k_1D^{(1)}_{\mathrm{cyl,aux}})e^{-\mathrm i n \phi^{(1)}_{\mathrm{aux}}}\mathrm{d}\phi^{(1)}_{\mathrm{aux}}\\
 &H_n^{(2)}(k_2\rho^{(2)}_{\mathrm{aux}})J'_n(k_2\rho_{\mathrm{cyl}})e^{-\mathrm i n \phi_{\mathrm{cyl}}}=\\ 
 &-\frac{1}{2\pi}\int_{-\pi}^{\pi}\frac{\rho_{\mathrm{cyl}}-\rho^{(2)}_{\mathrm{aux}}\cos(\phi_{\mathrm{cyl}}-\phi^{(2)}_{\mathrm{aux}})}{D^{(2)}_{\mathrm{cyl,aux}}}H_1^{(2)}(k_2D^{(2)}_{\mathrm{cyl,aux}})e^{-\mathrm i n \phi^{(2)}_{\mathrm{aux}}}\mathrm{d}\phi^{(2)}_{\mathrm{aux}}.
\end{align}
%

In this way, we obtain the algebraic system with respect to $I^s_n$ and $M^s_n$
\begin{align}
\label{eq:lin_sys_coeffs_1}
&I^s_nH_n^{(2)}(k_1\rho_{\mathrm{cyl}})+\frac{1}{\mathrm i Z_1}M^s_nH_n^{'(2)}(k_1\rho_{\mathrm {cyl}})=-\frac{I}{2\pi\rho_{\mathrm{cyl}}}H_n^{(2)}(k_1\rho_{\mathrm{fil}})\\
&I^s_nJ_n(k_2\rho_{\mathrm{cyl}})+\frac{1}{\mathrm i Z_2}M^s_nJ'_n(k_2\rho_{\mathrm {cyl}})=0
\label{eq:lin_sys_coeffs_2}
\end{align}
Importantly, the radii of the auxiliary surfaces (i.e., $\rho_\mathrm{aux}^{(1)}$ and $\rho_\mathrm{aux}^{(2)}$) have cancelled out and do not appear. Solving the system (\ref{eq:lin_sys_coeffs_1}) and (\ref{eq:lin_sys_coeffs_2}),  
 yields, via (\ref{eq:Fourier_series_1}) and  (\ref{eq:Fourier_series_2}), the explicit expressions of the current densities
\begin{align}
\nonumber
&J^s_z(\phi_{\mathrm{cyl}})=-\frac{I}{2\pi\rho_{\mathrm{cyl}}}\times 
\\
\label{eq:current_densities_ext_1}
&\sum_{n=-\infty}^{+\infty} \frac{Z_1H_n^{(2)}(k_1\rho_{\mathrm{fil}})J'_n(k_2\rho_{\mathrm{cyl}})}{Z_1H_n^{(2)}(k_1\rho_{\mathrm{cyl}})J'_n(k_2\rho_{\mathrm{cyl}})-Z_2H_n^{'(2)}(k_1\rho_{\mathrm{cyl}})J_n(k_2\rho_{\mathrm{cyl}})}e^{\mathrm i n \phi_{\mathrm{cyl}}}\\
\nonumber
&M^s_{\phi}(\phi_{\mathrm{cyl}})=\frac{I}{2\pi\rho_{\mathrm{cyl}}}\times 
\\
&\sum_{n=-\infty}^{+\infty}\frac{\mathrm i Z_1Z_2H_n^{(2)}(k_1\rho_{\mathrm{fil}})J_n(k_2\rho_{\mathrm{cyl}})}{Z_1H_n^{(2)}(k_1\rho_{\mathrm{cyl}})J'_n(k_2\rho_{\mathrm{cyl}})-Z_2H_n^{'(2)}(k_1\rho_{\mathrm{cyl}})J_n(k_2\rho_{\mathrm{cyl}})}e^{\mathrm i n \phi_{\mathrm{cyl}}}
\label{eq:current_densities_ext_2}
\end{align}

By the large-$n$ formulas of \ref{sec:appendix}, the $n$-th terms of the series in Eqs.~(\ref{eq:current_densities_ext_1}) and (\ref{eq:current_densities_ext_2}) behave, respectively, like 
\begin{align}
\label{eq:asymptotics_ext_1}
\frac{Z_1}{Z_1+Z_2\frac{k_2}{k_1}}\frac{\rho_{\mathrm{cyl}}^{|n|}}{\rho_{\mathrm{fil}}^{|n|}}e^{\mathrm i n \phi_{\mathrm{cyl}}},\\
-\mathrm i \frac{Z_1Z_2}{\frac{Z_1}{k_2}-\frac{Z_2}{k_1}}\frac{\rho^{|n|+1}_{\mathrm{cyl}}}{|n|\rho^{|n|}_{\mathrm{fil}}}e^{\mathrm i n \phi_{\mathrm{cyl}}}.
\label{eq:asymptotics_ext_2}
\end{align}
Thus, both series are convergent; this is true even if the auxiliary surfaces are placed beyond the extended regions mentioned in the previous section. In other words, the NFM integral equations---as opposed to those of MAS \cite{valagiannopoulos2012}---are always solvable.

Now, it is true that substitution of the expressions (\ref{eq:current_densities_ext_1}) and (\ref{eq:current_densities_ext_2}) for the surface current densities into the integrals (\ref{eq:R1_field}) for the fields do yield the exact solutions in (\ref{eq:exact_ext_TM_1}). As this is not immediately apparent (due to the complexity of the involved expressions), we provide a direct verification in \ref{sec:correct-fields}. 


By following the same steps for internal excitation ($\rho_{\mathrm{fil}}<\rho_{\mathrm{cyl}}$), we find
\begin{align}
\nonumber
&J_z^s(\phi_{\mathrm{cyl}})=-\frac{I}{2\pi\rho_{\mathrm{cyl}}}\times
\\ 
\label{eq:current_densities_int_1}
&\sum_{n=-\infty}^{+\infty} \frac{Z_2J_n(k_2\rho_{\mathrm{fil}})H^{'(2)}_n(k_1\rho_{\mathrm{cyl}})}{Z_2J_n(k_2\rho_{\mathrm{cyl}})H_n^{'(2)}(k_1\rho_{\mathrm {cyl}})-Z_1J'_n(k_2\rho_{\mathrm{cyl}})H^{(2)}_n(k_1\rho_{\mathrm{cyl}})}e^{\mathrm i n \phi_{\mathrm{cyl}}}\\
\nonumber
&M_{\phi}^s(\phi_{\mathrm{cyl}})=\frac{I}{2\pi\rho_{\mathrm{cyl}}}\times 
\\
&\sum_{n=-\infty}^{+\infty}\frac{\mathrm i Z_1Z_2J_n(k_2\rho_{\mathrm{fil}})H^{(2)}_n(k_1\rho_{\mathrm{cyl}})}{Z_2J_n(k_2\rho_{\mathrm{cyl}})H^{'(2)}_n(k_1\rho_{\mathrm{cyl}})-Z_1J'_n(k_2\rho_{\mathrm{cyl}})H^{(2)}_n(k_1\rho_{\mathrm{cyl}})}e^{\mathrm i n \phi_{\mathrm{cyl}}}.
\label{eq:current_densities_int_2}
\end{align}
The conclusions are similar. The series are always (irrespective of the positions of the auxiliary surfaces) convergent; and the current densities represented by these series always lead to the correct fields given in Section~\ref{sec:circular_problem}.

\section{Discrete Version of NFM-SEP}
\label{sec:Discrete_EIE}

\subsection{The Linear System and its Explicit Solution}
\label{subsec:The NFM Linear System and its Explicit Solution}

In the previous section, we considered the continuous version of NFM-SEP. In actual practice, we apply the discretized version described in Section~\ref{sec:NFM-SEP-discrete}. For the circular problem, we now formulate this system and obtain its closed-form solution. 
\begin{figure}[htb!]
    \centering
    \subfigure[]{\includegraphics[width=0.42\textwidth]{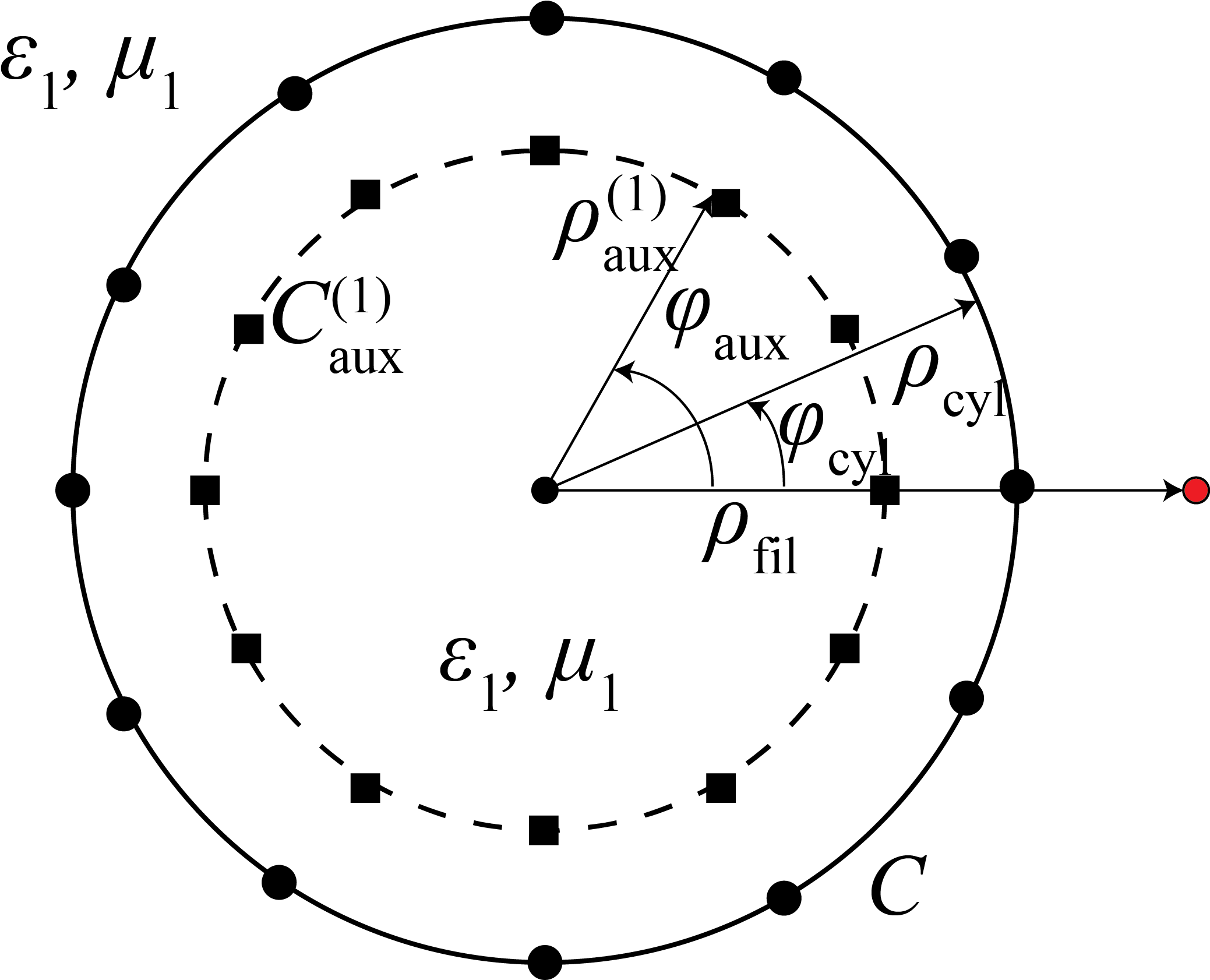}} 
    \subfigure[]{\includegraphics[width=0.415\textwidth]{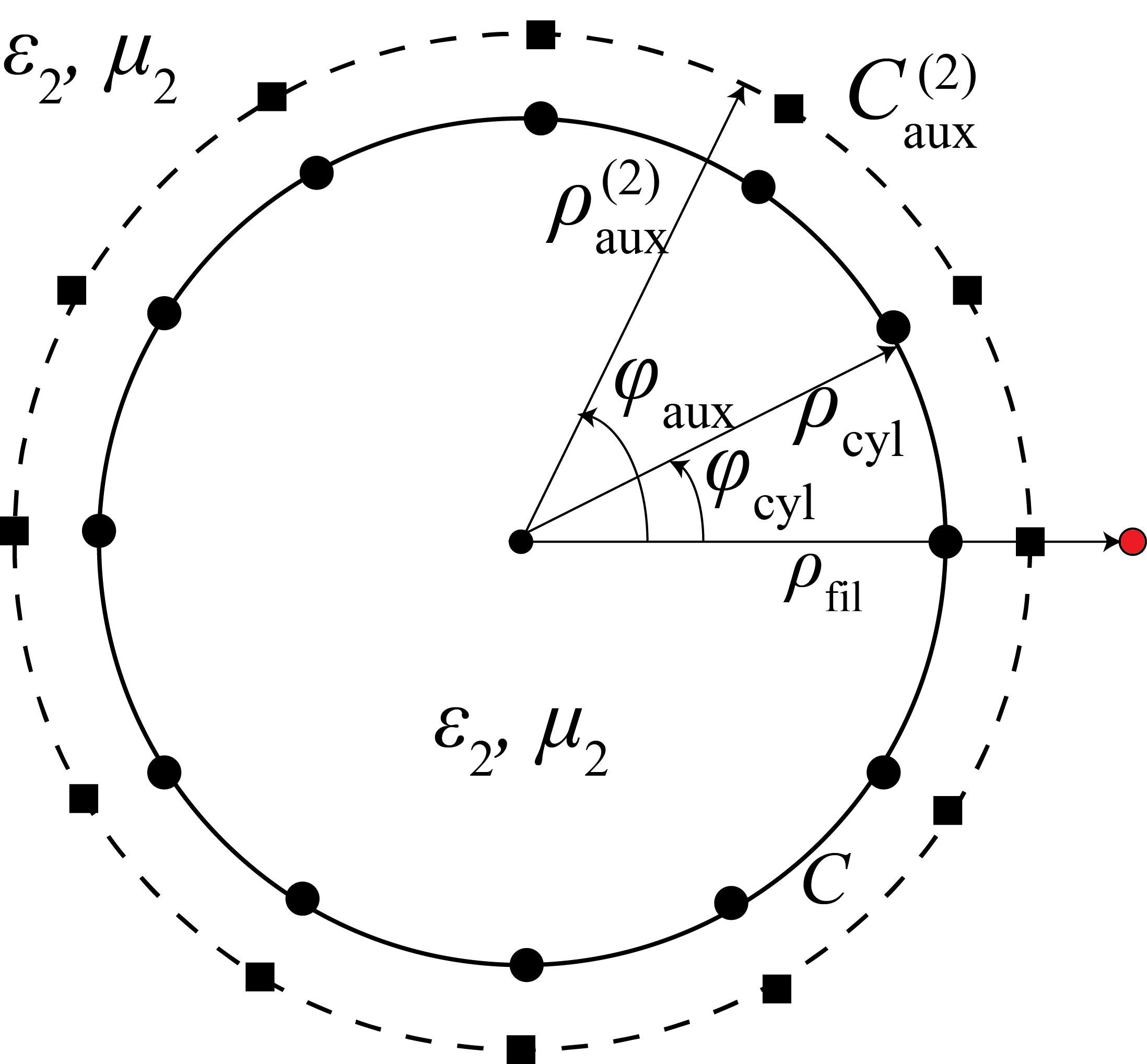}}
    \caption{Discretization of (a) the external and (b) the internal NFM equivalent problem for the circular dielectric cylinder.}
    \label{fig:Discrete_EIE}
\end{figure}
%


The continuous surface current densities of the previous section are replaced by discrete sources, as shown in Fig.~\ref{fig:Discrete_EIE}. 
Specifically, for external $\rm{TM}_z$ incidence, the electric current density $J_z^{s}(\phi_{\mathrm {cyl}})$ is replaced by $N$ $z$-directed electric current filaments with amplitudes $I_l$, and the magnetic current density $M_{\phi}^{s}(\phi_{\mathrm {cyl}})$ by $N$ $\phi$-directed magnetic current filaments with amplitudes $K_l$ ($l=0,1,\ldots,N-1$).
The amplitudes $I_l$ and $K_l$ are determined by satisfying the boundary conditions at $N$ collocation points on the polar angles $2p\pi/N$ of $C^{(1)}_{\mathrm {aux}}$ and $C^{(2)}_{\mathrm{aux}}$. For $p=0,1,\ldots,N-1$, this procedure leads to 
\begin{align}
\nonumber
&Z_1\sum_{l=0}^{N-1}I_lH_0^{(2)}(k_1b^{(1)}_{p-l})+\mathrm{i}\sum_{l=0}^{N-1}K_l\frac{\rho_{\mathrm{cyl}}-\rho^{(1)}_{\mathrm {aux}}\cos(\frac{2\pi(p-l)}{N})}{b^{(1)}_{p-l}}H_1^{(2)}(k_1b^{(1)}_{p-l})\\
&=-I Z_1 H_0^{(2)}(k_1d^{(1)}_p),
\label{eq:discrete_BC1}\\
&Z_2\sum_{l=0}^{N-1}I_lH_0^{(2)}(k_2b^{(2)}_{p-l})
+\mathrm{i}\sum_{l=0}^{N-1}K_l\frac{\rho_{\mathrm{cyl}}-\rho^{(2)}_{\mathrm {aux}}\cos(\frac{2\pi(p-l)}{N})}{b^{(2)}_{p-l}}H_1^{(2)}(k_2b^{(2)}_{p-l})=0,
\label{eq:discrete_BC2}
\end{align}
with
\begin{equation}
b^{(j)}_{p-l}=\sqrt{(\rho_{\mathrm {cyl}})^2+(\rho^{(j)}_{\mathrm{aux}})^2-2\rho_{\mathrm {cyl}}\rho^{(j)}_{\mathrm{aux}}\cos\left(\frac{2\pi(p-l)}{N}\right)}, \quad j=1,2
\end{equation}
the distance between the $p$-th collocation point on $C^{(j)}_{\mathrm{aux}}$ and the $l$-th source point on $C$, and
\begin{equation}
d^{(1)}_p=\sqrt{(\rho^{(1)}_{\mathrm {aux}})^2+(\rho_{\mathrm{fil}})^2-2\rho^{(1)}_{\mathrm {aux}}\rho_{\mathrm{fil}}\cos\left(\frac{2\pi p}{N}\right)}
\end{equation}
the distance from the excitation source to the $p$-th collocation point on $C^{(1)}_{\mathrm {aux}}$. 

Equations (\ref{eq:discrete_BC1}) and (\ref{eq:discrete_BC2})---which are  discretized versions of Eqs.~(\ref{eq:cont_eq_1}) and (\ref{eq:cont_eq_2})---are an $2N\times 2N$ linear system of equations with unknowns the complex amplitudes $I_l$ and $K_l$. In block form, the system is
\begin{equation}
\begin{bmatrix}
\mathbf Z_{11} & \mathbf Z_{12} \\ 
\mathbf Z_{21} & \mathbf Z_{22} 
\end{bmatrix}
\begin{bmatrix}
\mathbf I \\ 
\mathbf K
\end{bmatrix}=
\begin{bmatrix}
\mathbf V \\ 
\mathbf 0
\end{bmatrix},
\label{eq:NFM-block-matrix-system}
\end{equation}
where $\mathbf I$ and $\mathbf K$ are $N$-dimensional vectors with elements $I_l$ and $K_l$, respectively, while the $N$-dimensional vector $\mathbf V$ incorporates the contributions from the excitation source according to the right hand side of Eq.~(\ref{eq:discrete_BC1}). 
Importantly, each block $\mathbf Z_{11}$, $\mathbf Z_{12}$, $\mathbf Z_{21}$ and $\mathbf Z_{22}$ is a circulant matrix, namely,  a square matrix in which all rows are composed of the same elements and each row is rotated one element to the right relative to the preceding row \cite{Davis1994, Gray2002}. Hence the matrix in  (\ref{eq:NFM-block-matrix-system}) is a $2\times 2$ block matrix composed of circulant blocks of order $N$. 
Block matrices with circulant blocks were analyzed in
\cite{Tsitsas2010}, where efficient inversion algorithms were established.


We solve (\ref{eq:NFM-block-matrix-system}) by utilizing the discrete Fourier transform (DFT)  (similar to the  procedure in \cite{valagiannopoulos2012}).  We consider the following DFT pairs of $I_l$ and $K_l$,
%
%
\begin{equation}
I_l=\sum_{m=0}^{N-1}I^{(m)}e^{\mathrm i 2\pi l m/N},\,\,
K_l=\sum_{m=0}^{N-1}K^{(m)}e^{\mathrm i 2\pi l m/N},~l=0,1,\ldots,N-1
\label{eq:solution_CS}
\end{equation}
\begin{equation}
I^{(m)}=\frac{1}{N}\sum_{p=0}^{N-1}I_pe^{-\mathrm i 2\pi m p/N},\,\,
K^{(m)}=\frac{1}{N}\sum_{p=0}^{N-1}K_pe^{-\mathrm i 2\pi mp/ N},~m=0,1,\ldots,N-1.
\label{eq:current_fourier_coefs}
\end{equation} 

The DFT of any convolution sum $\sum_{l=0}^{N-1}\alpha_{p-l}\beta_l$ with $\alpha_l=\alpha_{-l}=\alpha_{N-l}$ equals $N\alpha^{(m)}\beta^{(m)}$. Thus, application of the DFT to (\ref{eq:discrete_BC1}) and (\ref{eq:discrete_BC2}) gives the following $2\times 2$ system 
\begin{align}
\label{eq:discrete_LS_1}    
&B^{(m)}_1 I^{(m)}+\frac{\mathrm i}{Z_1}B_2^{(m)}K^{(m)}=\frac{I}{N} D^{(m)}\\
&B_3^{(m)}I^{(m)}+\frac{\mathrm i}{ Z_2}B_4^{(m)}K^{(m)}=0,
\label{eq:discrete_LS_2}
\end{align}
for all $m=0,1,\ldots,N-1$, with
\begin{align}
\label{eq:B1m}
B_1^{(m)}&=\frac{1}{N}\sum_{p=0}^{N-1}H_0^{(2)}(k_1b^{(1)}_{p})e^{-\mathrm i 2\pi m p/N},\\
\label{eq:B2m}
B_2^{(m)}&=\frac{1}{N}\sum_{p=0}^{N-1}\frac{\rho_{\mathrm{cyl}}-\rho^{(1)}_{\mathrm{aux}}\cos(\frac{2\pi p}{N})}{b^{(1)}_{p}}H_1^{(2)}(k_1b^{(1)}_{p})e^{-\mathrm i 2\pi m p/N},\\
\label{eq:B3m}
B_3^{(m)}&=\frac{1}{N}\sum_{p=0}^{N-1}H_0^{(2)}(k_2b^{(2)}_{p})e^{-\mathrm i 2\pi m p/N},\\
\label{eq:B4m}
B_4^{(m)}&=\frac{1}{N}\sum_{p=0}^{N-1}\frac{\rho_{\mathrm{cyl}}-\rho^{(2)}_{\mathrm{aux}}\cos(\frac{2\pi p}{N})}{{b^{(2)}_{p}}}H_1^{(2)}(k_2b^{(2)}_{p})e^{-\mathrm i 2\pi m p/N},\\
\label{eq:Dm}
D^{(m)}&=-\frac{1}{N}\sum_{p=0}^{N-1}H_0^{(2)}(k_1d_{p}^{(1)})e^{-\mathrm i 2\pi m p/N}.
\end{align}

Then, by solving (\ref{eq:discrete_LS_1}) and (\ref{eq:discrete_LS_2}), we obtain
\begin{align}
I^{(m)}&=\frac{I}{N}\frac{Z_1D^{(m)}B_4^{(m)}}{Z_1B_1^{(m)}B_4^{(m)}-Z_2B_2^{(m)}B_3^{(m)}}\\
K^{(m)}&=\frac{I}{N}\frac{\mathrm i Z_2Z_1D^{(m)}B_3^{(m)}}{Z_1B_1^{(m)}B_4^{(m)}-Z_2B_2^{(m)}B_3^{(m)}},
\label{eq:DFT_currents}
\end{align}
which, when combined with (\ref{eq:solution_CS}), provide the explicit, closed-form solution of the system (\ref{eq:discrete_BC1}) and (\ref{eq:discrete_BC2}), for any finite $N$. We now study  the large-$N$ behavior of this solution.

\subsection{Discrete Currents and Fields they Generate: Large-$N$ behavior}

We specifically ask whether the two limits
\begin{equation}
J^s_{\rm{limit}}(\phi_{\mathrm{cyl}})=\lim_{N \to \infty}\frac{NI_l}{2\pi\rho_{\mathrm{cyl}}}\,\,\,\mbox{and}\,\,\,M^s_{\rm{limit}}(\phi_{\mathrm{cyl}})=\lim_{N \to \infty}\frac{NK_l}{2\pi\rho_{\mathrm{cyl}}}
\label{eq:normalized_currents}
\end{equation}
exist; and, if they do exist, whether they coincide with the surface current densities (\ref{eq:current_densities_ext_1}) and (\ref{eq:current_densities_ext_2}). 
The question thus posed is a natural one because \cite{valagiannopoulos2012} shows that, in the case of MAS, the analogous limits do not always exist.

We proceed  as follows: Replace the Hankel functions in (\ref{eq:B1m})--(\ref{eq:Dm}) by
the expressions given in (\ref{eq:addition_theorem})-(\ref{eq:addition_theorem_der_2}); interchange the order of summation in the resulting expressions; and then use the identity
\begin{equation}
\sum_{p=0}^{N-1}e^{\mathrm i 2\pi p(n-m)/N}=
\begin{cases}
N,~~\rm{if~n-m=multiple~of~N}\\
0,~~\rm{otherwise}
\end{cases}
\end{equation}
This procedure results in the exact expressions (for $m=0,\ldots,N-1$)
\begin{align}
\label{eq:q-sums-1}
D^{(m)}&=-\sum_{q=-\infty}^{+\infty}J_{qN+m}(k_1\rho^{(1)}_{\mathrm{aux}})H^{(2)}_{qN+m}(k_1\rho_{\mathrm{fil}}),\\
\label{eq:q-sums-2}
B^{(m)}_1&=\sum_{q=-\infty}^{+\infty}J_{qN+m}(k_1\rho^{(1)}_{\mathrm{aux}})H^{(2)}_{qN+m}(k_1\rho_{\mathrm{cyl}}),\\
\label{eq:q-sums-3}
B^{(m)}_2&=-\sum_{q=-\infty}^{+\infty}J_{qN+m}(k_1\rho^{(1)}_{\mathrm{aux}})H'^{(2)}_{qN+m}(k_1\rho_{\mathrm{cyl}}),\\
\label{eq:q-sums-4}
B^{(m)}_3&=\sum_{q=-\infty}^{+\infty}J_{qN+m}(k_2\rho_{\mathrm{cyl}})H^{(2)}_{qN+m}(k_2\rho^{(2)}_{\mathrm{aux}}),\\
\label{eq:q-sums-5}
B^{(m)}_4&=-\sum_{q=-\infty}^{+\infty}J'_{qN+m}(k_2\rho_{\mathrm{cyl}})H^{(2)}_{qN+m}(k_2\rho^{(2)}_{\mathrm{aux}}).
\end{align}

Now, assume for simplicity that $N$ is odd, so that (\ref{eq:solution_CS}) is written as
\begin{align}
\label{eq:solution_CS_odd_1}
I_l&=I^{(0)}+\sum_{m=1}^{(N-1)/2}I^{(m)}\cos\left(\frac{2\pi lm}{N}\right),~l=0,\ldots,N-1,\\
K_l&=K^{(0)}+\sum_{m=1}^{(N-1)/2}K^{(m)}\cos\left(\frac{2\pi lm}{N}\right),~l=0,\ldots,N-1.
\label{eq:solution_CS_odd_2}
\end{align}

Then, for large $N$ and $0 \le m \le (N-1)/2$, by (\ref{eq:Bessel})-(\ref{eq:Hankel_der}), it is seen that the $q=0$ terms in (\ref{eq:q-sums-1})-(\ref{eq:q-sums-5}) are dominant. Thus, as $N\to\infty$, we get
\begin{align}
\label{eq:large_forms_1}
I^{(m)}\sim -\frac{I}{N}\frac{Z_1H_m^{(2)}(k_1\rho_{\mathrm{fil}})J'_m(k_2\rho_{\mathrm{cyl}})}{Z_1H_m^{(2)}(k_1\rho_{\mathrm{cyl}})J'_m(k_2\rho_{\mathrm{cyl}})-Z_2H_m^{'(2)}(k_1\rho_{\mathrm{cyl}})J_m(k_2\rho_{\mathrm{cyl}})},\\
K^{(m)}\sim \frac{I}{N}\frac{\mathrm i Z_1Z_2H_m^{(2)}(k_1\rho_{\mathrm{fil}})J_m(k_2\rho_{\mathrm{cyl}})}{Z_1H_m^{(2)}(k_1\rho_{\mathrm{cyl}})J'_m(k_2\rho_{\mathrm{cyl}})-Z_2H_m^{'(2)}(k_1\rho_{\mathrm{cyl}})J_m(k_2\rho_{\mathrm{cyl}})},
\label{eq:large_forms_2}
\end{align}
in which $\rho_\mathrm{aux}^{(1)}$ and $\rho_\mathrm{aux}^{(2)}$ have cancelled out. Upon substituting these expressions into 
(\ref{eq:solution_CS_odd_1}) and (\ref{eq:solution_CS_odd_2}) and setting $2\pi l/N=\phi_\mathrm{cyl}$, we see that the limits defined in (\ref{eq:normalized_currents}) are finite and given by the series
\begin{align}
\nonumber
&J^s_{\rm{limit}}(\phi_{\mathrm{cyl}})=-\frac{I}{2\pi\rho_{\mathrm{cyl}}}\times\\
\nonumber
&\biggl[\frac{Z_1H_0^{(2)}(k_1\rho_{\mathrm{fil}})J'_0(k_2\rho_{\mathrm{cyl}})}{Z_1H_0^{(2)}(k_1\rho_{\mathrm{cyl}})J'_0(k_2\rho_{\mathrm{cyl}})-Z_2H_0^{'(2)}(k_1\rho_{\mathrm{cyl}})J_0(k_2\rho_{\mathrm{cyl}})}+\\ 
&2\sum_{m=1}^{+\infty} \frac{Z_1H_m^{(2)}(k_1\rho_{\mathrm{fil}})J'_m(k_2\rho_{\mathrm{cyl}})}{Z_1H_m^{(2)}(k_1\rho_{\mathrm{cyl}})J'_m(k_2\rho_{\mathrm{cyl}})-Z_2H_m^{'(2)}(k_1\rho_{\mathrm{cyl}})J_m(k_2\rho_{\mathrm{cyl}})}\cos(  m \phi_{\mathrm{cyl}}) \biggr],\\
\nonumber
&M^s_{\rm{limit}}(\phi_{\mathrm{cyl}})=\frac{I}{2\pi\rho_{\mathrm{cyl}}}\times\\
\nonumber
&\biggl[\frac{\mathrm i Z_1Z_2H_0^{(2)}(k_1\rho_{\mathrm{fil}})J_0(k_2\rho_{\mathrm{cyl}})}{Z_1H_0^{(2)}(k_1\rho_{\mathrm{cyl}})J'_0(k_2\rho_{\mathrm{cyl}})-Z_2H_0^{'(2)}(k_1\rho_{\mathrm{cyl}})J_0(k_2\rho_{\mathrm{cyl}})}+\\ 
&2\sum_{m=1}^{+\infty} \frac{\mathrm i Z_1Z_2H_m^{(2)}(k_1\rho_{\mathrm{fil}})J_m(k_2\rho_{\mathrm{cyl}})}{Z_1H_m^{(2)}(k_1\rho_{\mathrm{cyl}})J'_m(k_2\rho_{\mathrm{cyl}})-Z_2H_m'^{(2)}(k_1\rho_{\mathrm{cyl}})J_m(k_2\rho_{\mathrm{cyl}})}\cos( m \phi_{\mathrm{cyl}}) \biggr],
\end{align}
which are Fourier series whose coefficients are even in $m$.
Therefore,  $J^s_{\rm{limit}}(\phi_{\mathrm{cyl}})$ and $M^s_{\rm{limit}}(\phi_{\mathrm{cyl}})$ are equal to the $J_z^s(\phi_{\mathrm{cyl}})$ and $M_\phi^s(\phi_{\mathrm{cyl}})$ given by (\ref{eq:current_densities_ext_1}) and (\ref{eq:current_densities_ext_2}), respectively. 
Hence, it follows that the large-$N$ limits of the fields generated by $I_l$ and $K_l$ exist and equal the fields generated by $J_z^s(\phi_{\mathrm{cyl}})$ and $M_\phi^s(\phi_{\mathrm{cyl}})$. Recall that these were shown in \ref{sec:correct-fields} to be the true fields.

Repeating the calculations for internal $\rm{TM_z}$ excitation, we reached the same conclusions.

\section{NFM-SEP vs. MAS}
\label{nfm-sep-vs-mas}

We summarize our hitherto results by comparing to corresponding ones which can be found in \cite{valagiannopoulos2012} and which pertain to MAS. Numerical methods are often compared between themselves, but the comparison we now perform (between NFM-SEP and MAS) is more appropriate than usual; this is so because the matrices of the two methods are, apart from trivial differences  \footnote{The trivial differences follow by comparing (\ref{eq:discrete_BC1}) and (\ref{eq:discrete_BC2}) to (11) of \cite{valagiannopoulos2012}: We obtain the MAS matrix as follows. Multiply the top two blocks (i.e., the first $N$ rows) of our NFM-SEP matrix by $-k_1/4$; multiply the bottom two blocks by $k_2/4$; block-transpose. The multiplications of the blocks of the NFM-SEP matrix by the above coefficients correspond to respective multiplications on both sides of (\ref{eq:discrete_BC1}) and (\ref{eq:discrete_BC2}).}, the same. Only the right-hand-side vectors present nontrivial differences. Importantly, we note that this coincidence of the NFM-SEP and MAS matrices continues to hold also for non-circular shapes. 

\begin{itemize}
    \item{Continuous integral equations:} The NFM-SEP continuous integral equations are always solvable; but those of MAS are nonsolvable if the interior surface has radius smaller than a certain threshold value, or if the exterior surface has radius larger than an other threshold value; all thresholds depend on the type of excitation (external or internal) and are given in Table~\ref{table:divergence-MAS-currents}.
    \item{Solution to $2N\times 2N$ system:} In the case of NFM-SEP, the aforesaid solution (discrete currents), when appropriately normalized, always converge as $N\to\infty$. By contrast, the interior (exterior) MAS currents diverge if the interior (exterior) surface is too small (too large).  Since the $2N\times 2N$ system can be viewed as a discretization of the system of continuous integral equations, the divergence of the normalized solution to the former system is a simple corollary of the nonsolvability of the latter system.
    \item{Oscillations:} In MAS, the aforementioned divergence of the normalized currents manifests itself as rapid, unnatural, exponentially large oscillations in the corresponding MAS currents \cite{valagiannopoulos2012}. On the other hand, NFM-SEP does not give rise to such oscillations.
    \item{Field convergence:} In both NFM-SEP and MAS, the fields generated by the discrete currents always converge to the true field (regardless of the positions of the auxiliary surfaces and of the occurrence/nonoccurrence of oscillations).
    \item{Numerical issues:} Field convergence is a theoretical result, which might not hold true when the hardware and/or software are imperfect, which is always the case in practice. This is why the lack of oscillations in NFM-SEP is an advantage of NFM-SEP (when compared to MAS), in accordance with the principle that algorithms involving large intermediate results are generally unstable \cite{Fikioris2018convergent,higham} 
    \item{Singularities of analytic continuation:} The extended convergence regions we showed in Table~\ref{table:convergence} are in fact regions into which the physical fields can be analytically continued. As discussed in \cite{valagiannopoulos2012}, the singularities of these analytic continuations play an important role to the convergence/divergence of MAS. As shown conclusively in the present paper, the aforementioned singularities are \textit{not} relevant to the convergence/divergence of NFM-SEP. This finding is a generalization of a similar finding for perfectly conducting (PEC) scatterers \cite{fikioris2011}. In fact, for the PEC case the issue had caused some controversy: see \cite{fikioris2011,fikioris2013,Fikioris2018convergent,EreminSkobolev2011} and the relevant references therein.
\end{itemize}

\section{Numerical Results}
\label{sec:Numerical}

In this section, we show representative numerical results illustrating the main conclusions of the paper, both for circular as well as for elliptical cylinders. 
Throughout, 
the relative permittivity of the $R_2$ medium is taken as 
$\varepsilon_{r2}=4.2$. 
As in the theoretical part, we consider only $\rm{TM_z}$ excitation (both internal and external); for $\rm{TE_z}$ excitation the conclusions were similar.

\subsection{Circular Geometry}
\label{sec:Numerical-circular}

\begin{figure}[H]
    \centering
    \subfigure[]{\includegraphics[width=0.47\textwidth]{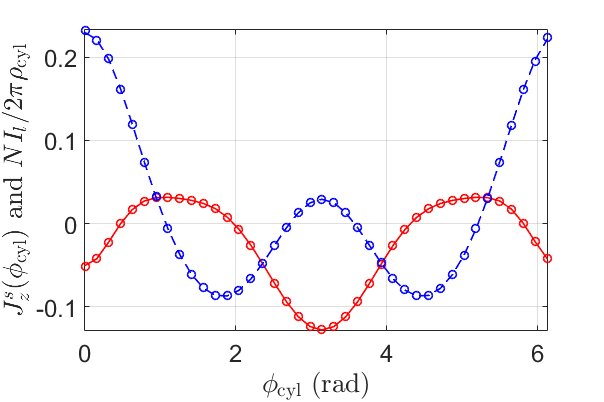}} 
    \subfigure[]{\includegraphics[width=0.47\textwidth]{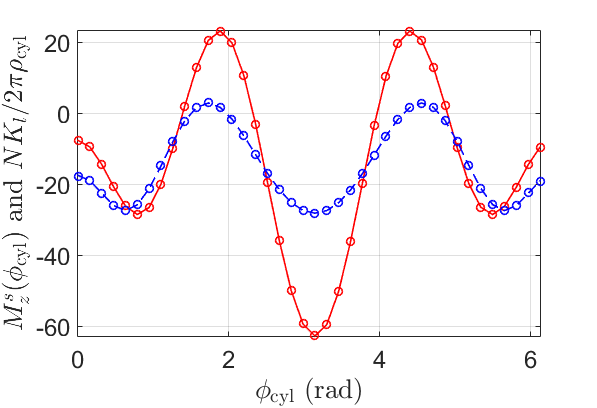}}
    \centering
    \subfigure[]{\includegraphics[width=0.47\textwidth]{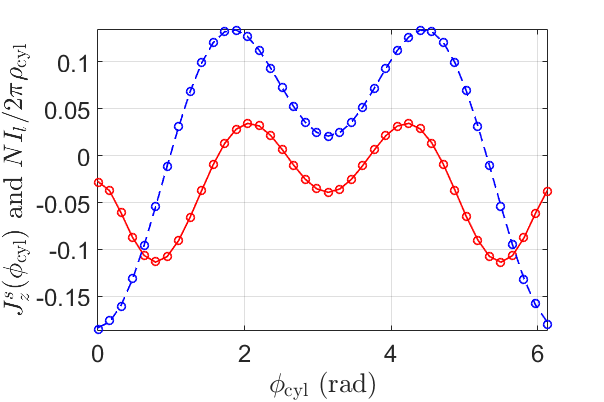}} 
    \subfigure[]{\includegraphics[width=0.47\textwidth]{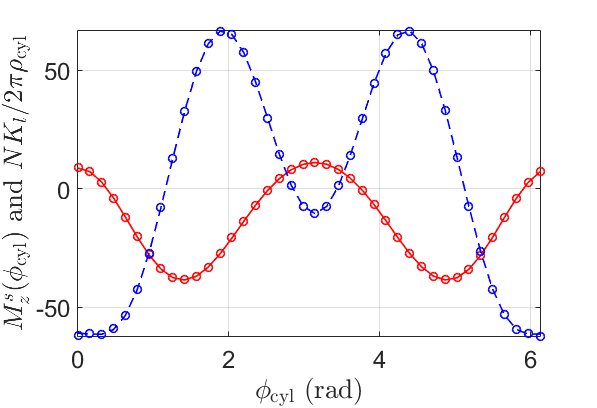}}
    \subfigure[]{\includegraphics[width=0.47\textwidth]{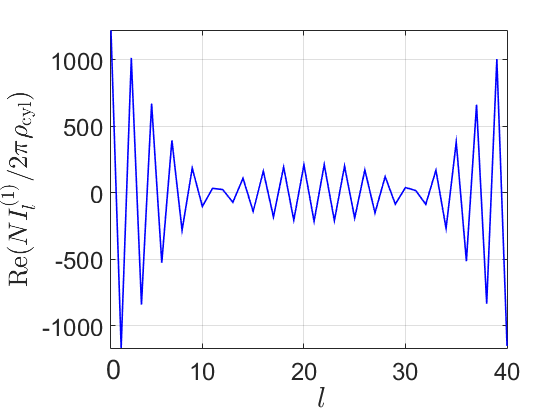}} 
    \subfigure[]{\includegraphics[width=0.47\textwidth]{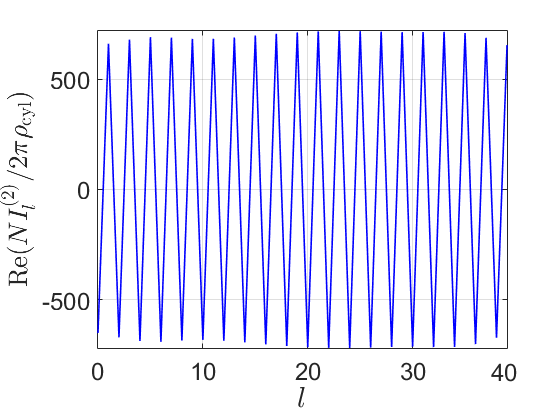}}
    \caption{Real (red solid lines) and imaginary parts (blue dashed lines) of (a) continuous electric current density $J_z^s(\phi_{\mathrm{cyl}})$ and (b) continuous magnetic current density $M_\phi^s(\phi_{\mathrm{cyl}})$ together with real and imaginary parts (circles) of normalized NFM currents $NI_l/2\pi\rho_{\mathrm{cyl}}$ and $NK_l/2\pi\rho_{\mathrm{cyl}}$ for external $\rm{TM_z}$  excitation. (c) and (d) are as (a) and (b) but for internal $\rm{TM_z}$ excitation.  For (a)-(d): $k_1\rho^{(1)}_{\mathrm {aux}}=1.5$, $k_1\rho^{(2)}_{\mathrm {aux}}=2.5$ (or $k_1\rho^{(1)}_{\mathrm {aux}}=0.5$, $k_1\rho^{(2)}_{\mathrm {aux}}=10$) and $N=40$. (e) and (f) show the oscillations of the MAS currents of $C^{(1)}_{\mathrm{aux}}$ for internal and  $C^{(2)}_{\mathrm{aux}}$ for external excitation with $k_1\rho^{(1)}_{\mathrm {aux}}=0.5$, $k_1\rho^{(2)}_{\mathrm {aux}}=10$ and $N=40$.}
    \label{fig:dens_vs_norm_TM}
\end{figure}

\begin{figure}[H]
    \centering
    \subfigure[]{\includegraphics[width=0.47\textwidth]{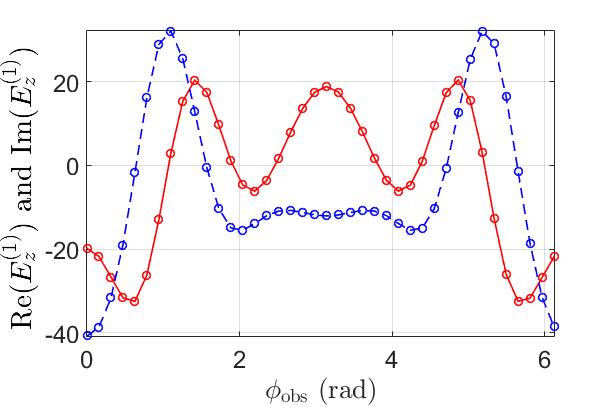}} 
    \subfigure[]{\includegraphics[width=0.47\textwidth]{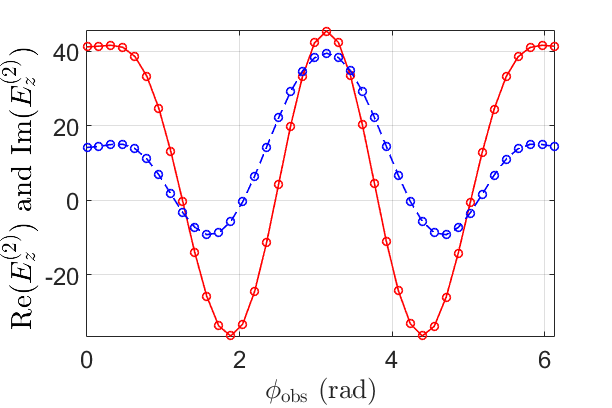}}
    \subfigure[]{\includegraphics[width=0.47\textwidth]{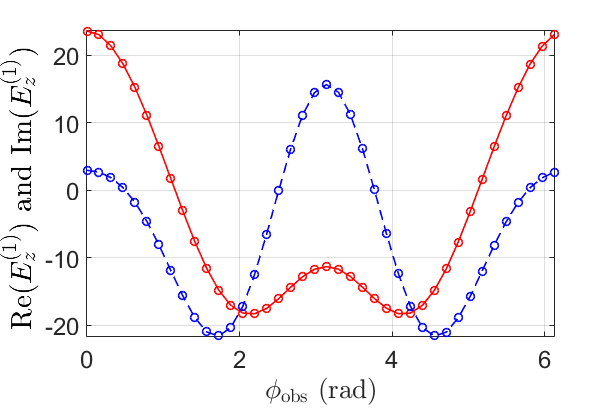}}
    \subfigure[]{\includegraphics[width=0.47\textwidth]{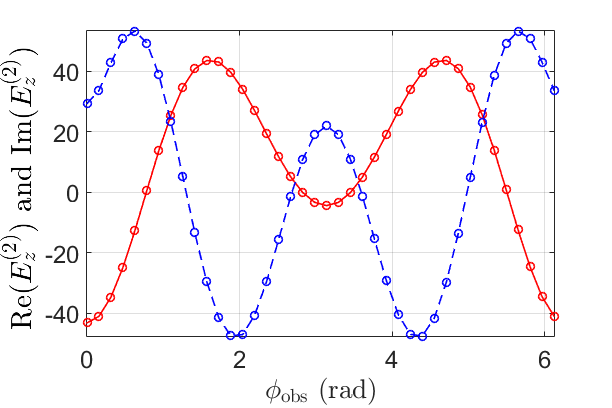}}
    \subfigure[]{\includegraphics[width=0.47\textwidth]{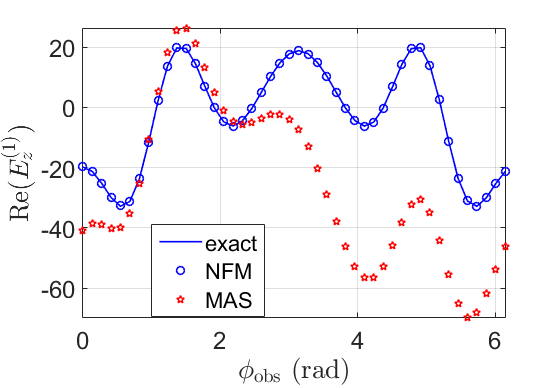}}
    \subfigure[]{\includegraphics[width=0.47\textwidth]{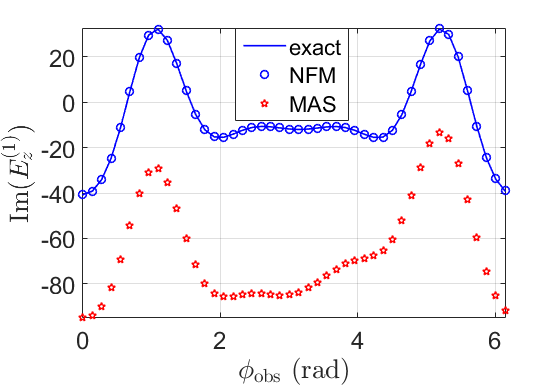}}
    \caption{Real (red solid lines) and imaginary parts (blue dashed lines) of the exact electric fields (a) $E_z^{(1)}(\phi_{\mathrm{obs}})$ and (b) $E_z^{(2)}(\phi_{\mathrm{obs}})$ together with real and imaginary parts of (a) $E_{z,\mathrm{NFM}}^{(1)}(\phi_{\mathrm{obs}})$ and (b) $E_{z,\mathrm{NFM}}^{(2)}(\phi_{\mathrm{obs}})$ (circles) for external $\rm{TM_z}$ excitation with $N=40$, $k_1\rho^{(1)}_{\mathrm {aux}}=0.5$, $k_1\rho^{(2)}_{\mathrm {aux}}=10$, $k_1\rho^{(1)}_{\mathrm {obs}}=10$, and $k_1\rho^{(2)}_{\mathrm {obs}}=1$. (c) and (d) are as (a) and (b) but for internal $\rm{TM_z}$ excitation. (e) and (f) illustrate the negative effect of oscillations (for $N=46$) of the MAS currents of $C^{(1)}_{\mathrm{aux}}$ on the computations of the real and imaginary parts of $E_z^{(1)}$ for external excitation.}
    \label{fig:exact_EIE_el_fields}
\end{figure}

Figures \ref{fig:dens_vs_norm_TM} and \ref{fig:exact_EIE_el_fields} concern the circular problem with $N=40$, $k_1\rho_{\mathrm {cyl}}=2$, $k_1\rho_{\mathrm{fil}}=4$ for external excitation, and $k_1\rho_{\mathrm{fil}}=1$ for internal excitation. Thus, $k_1\rcri=1$ and $k_1\rcri=4$ for external and internal excitation, respectively.
Figures \ref{fig:dens_vs_norm_TM}(a)-(d) compare the continuous surface current densities $J^s_z(\phi_{\mathrm{cyl}})$ and $M^s_\phi(\phi_{\mathrm{cyl}})$ of Section~\ref{sec:continuous_EIE} with the normalized currents $NI_l/(2\pi\rcyl)$ and $NK_l/(2\pi\rcyl)$ of Section~\ref{sec:Discrete_EIE},  for $k_1\rho^{(1)}_{\mathrm {aux}}=1.5$ and $k_1\rho^{(2)}_{\mathrm {aux}}=2.5$. Perfect coincidence is observed. 
These plots remain unaltered, when we change the radii of the auxiliary surfaces to $k_1\rho^{(1)}_{\mathrm {aux}}=0.5$ and $k_1\rho^{(2)}_{\mathrm {aux}}=10$. 
This verifies numerically that the NFM-SEP currents always converge to the expected values. 
In stark contrast, the real parts of the corresponding MAS currents on $C^{(1)}_{\mathrm{aux}}$ and $C^{(2)}_{\mathrm{aux}}$ (denoted, respectively, as $NI^{(1)}_{l}/2\pi \rho_{\mathrm{cyl}}$ and $NI^{(2)}_{l}/2\pi \rho_{\mathrm{cyl}}$) oscillate rapidly, as shown in Figs.~\ref{fig:dens_vs_norm_TM} (e)-(f). It is worth stressing that these oscillations are hardware/software \emph{independent} \cite{fikioris2006, Fikioris2018convergent}. 
The resulting NFM-SEP fields are compared to the exact fields in Figs.~\ref{fig:exact_EIE_el_fields}(a)-(d). 
Once again the results are coincident; they remain the same if we change the positions of $C^{(1)}_{\mathrm{aux}}$ and $C^{(2)}_{\mathrm{aux}}$ to $k_1\rho^{(1)}_{\mathrm {aux}}=0.5$ and $k_1\rho^{(2)}_{\mathrm {aux}}=10$. 
The MAS fields, \textit{despite the oscillations} shown in Figs.~\ref{fig:dens_vs_norm_TM}(e)-(f), also coincide for $N=40$ with the exact fields; this phenomenon is explained in detail in \cite{valagiannopoulos2012}. 

If $N$ is increased to $N=46$, the NFM-SEP fields remain correct. For the MAS solution, on the other hand, this small increase in $N$ causes a large increase in the amplitude of the oscillating currents $I^{(2)}_l$,  resulting in corrupted field calculations; see Figs.~\ref{fig:exact_EIE_el_fields}(e)-(f). 
We point out that such corrupted results are hardware/software \emph{dependent}. This illustrates the numerical issues discussed in Section~\ref{nfm-sep-vs-mas}, and  showcases the advantage of NFM-SEP over MAS.

\subsection{Elliptical Geometry}
\label{sec:Numerical-ellipse}

\begin{figure}[H]
    \centering
    \subfigure[]{\includegraphics[width=0.47\textwidth]{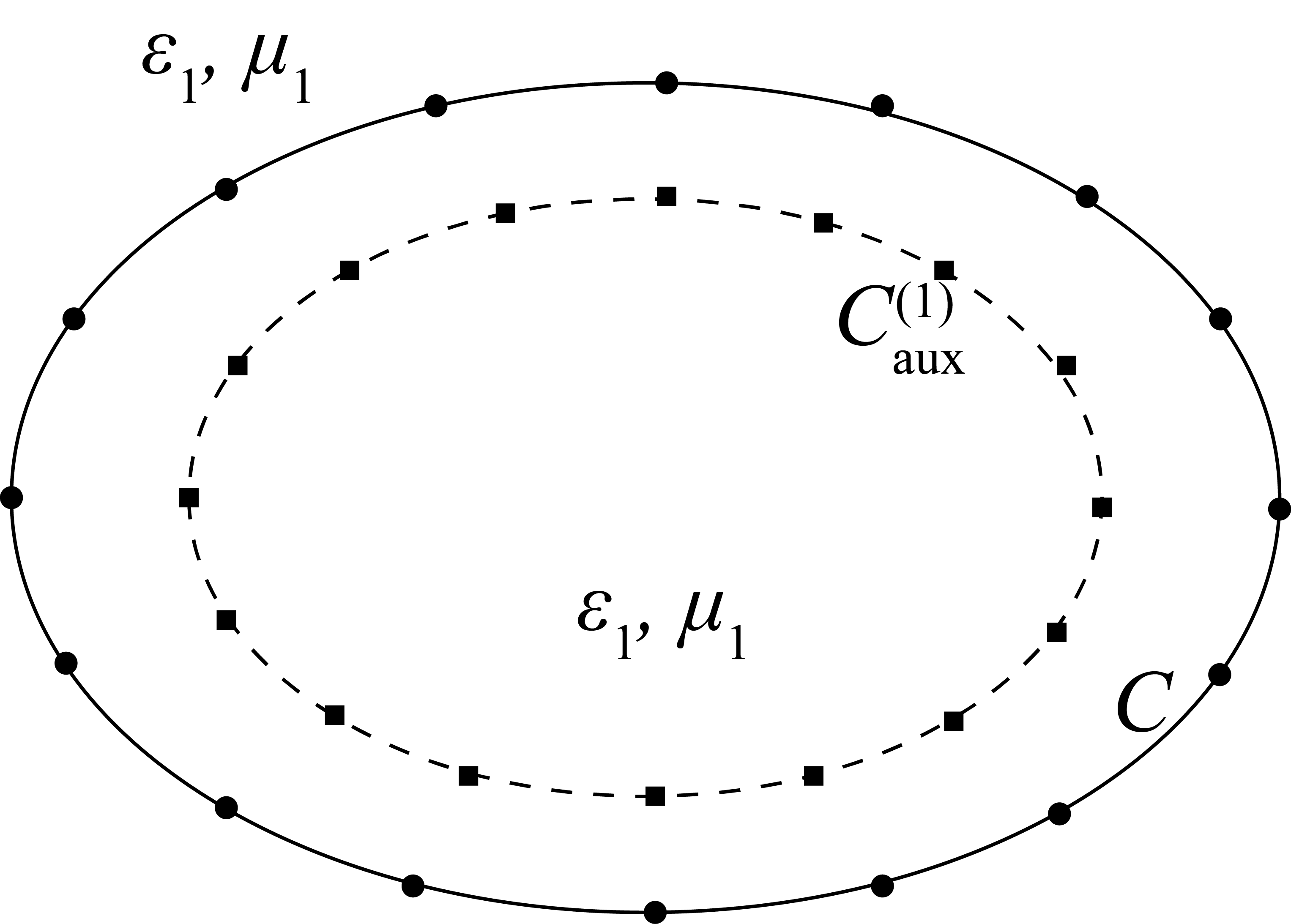}} 
    \subfigure[]{\includegraphics[width=0.50\textwidth]{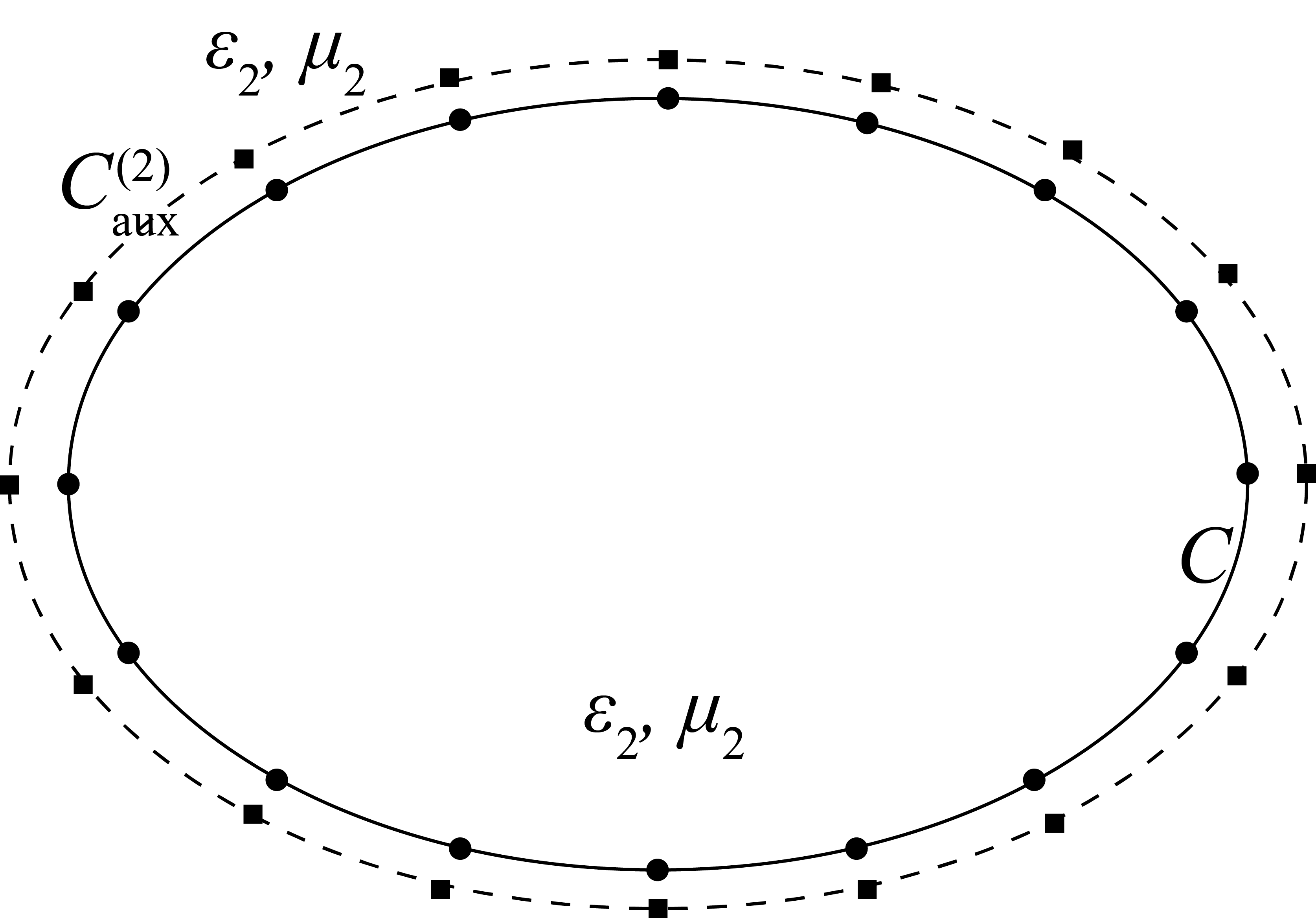}}
    \caption{Discretization of the (a) external and (b) internal NFM equivalent problem for the elliptical dielectric cylinder. Black dots represent the collocated electric and magnetic discrete sources, and black squares represent the collocation points.}
    \label{fig:ellipse_EIE}
\end{figure}

\begin{figure}[H]
    \centering
    \subfigure[]{\includegraphics[width=0.47\textwidth]{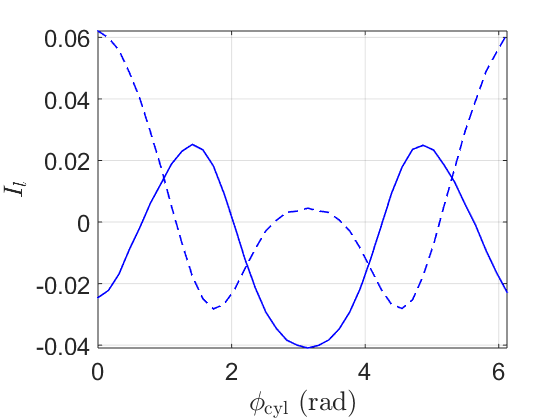}} 
    \subfigure[]{\includegraphics[width=0.47\textwidth]{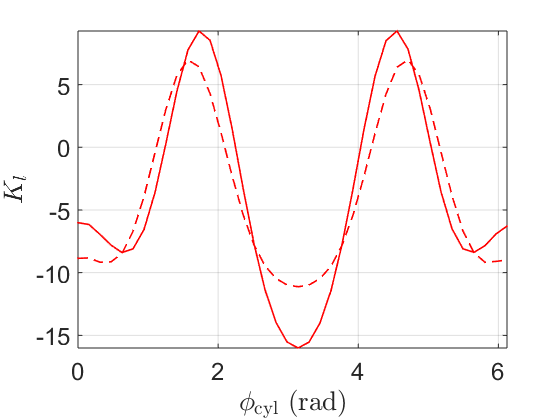}}
    \centering
    \subfigure[]{\includegraphics[width=0.47\textwidth]{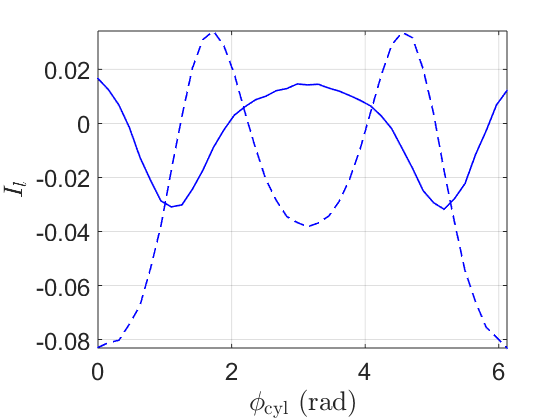}} 
    \subfigure[]{\includegraphics[width=0.47\textwidth]{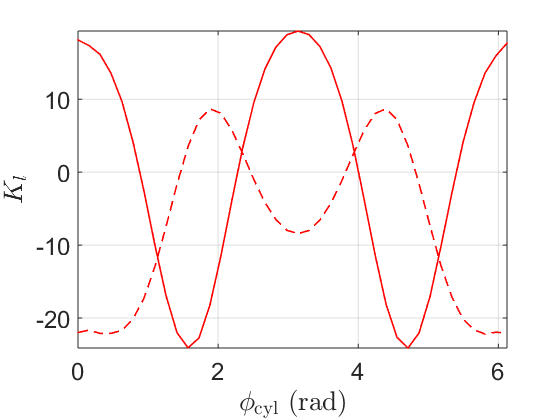}}
    \subfigure[]{\includegraphics[width=0.47\textwidth]{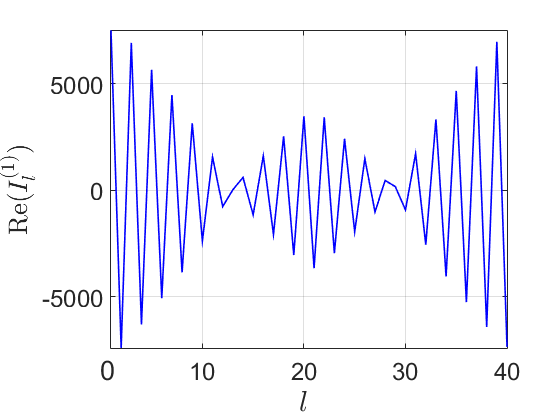}} 
    \subfigure[]{\includegraphics[width=0.47\textwidth]{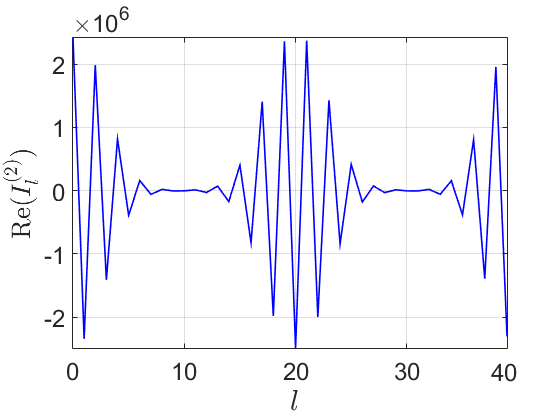}}
    \caption{Real (solid lines) and imaginary parts (dashed lines) of the NFM (a) electric current amplitudes $I_l$ and (b) magnetic current amplitudes $K_l$ for external $\rm {TM_z}$ excitation of the elliptical cylinder for $N=40$, $\sigma^{(1)}_{\mathrm {aux}}=0.33$ and $\sigma^{(2)}_{\mathrm {aux}}=5$. (c) and (d) are as (a) and (b) but  for internal $\rm{TM_z}$ excitation. (e) and (f) illustrate, respectively, the oscillations of the MAS currents of $C^{(1)}_{\mathrm{aux}}$ for external and $C^{(2)}_{\mathrm{aux}}$ for internal excitation.}
    \label{fig:dens_TM_ell}
\end{figure}

\begin{figure}[H]
    \centering
    \subfigure[]{\includegraphics[width=0.47\textwidth]{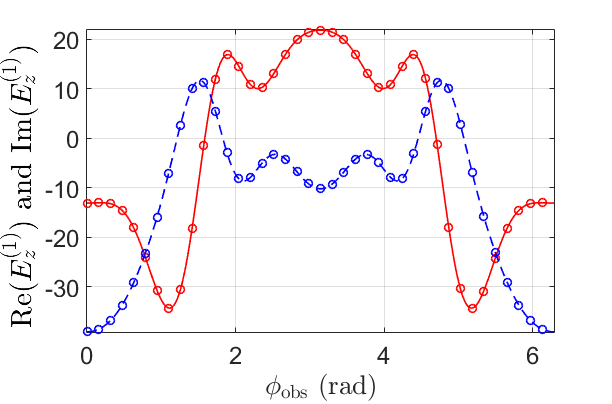}} 
    \subfigure[]{\includegraphics[width=0.47\textwidth]{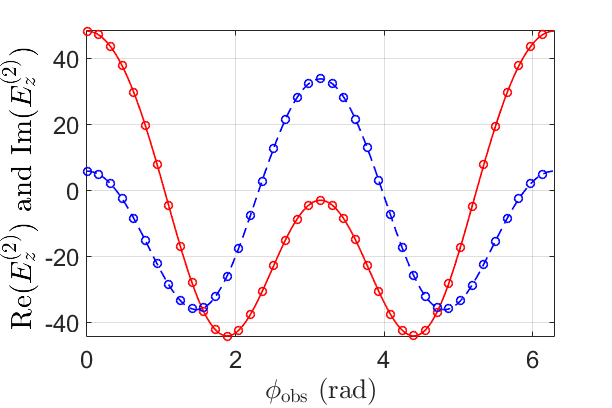}}
    \subfigure[]{\includegraphics[width=0.47\textwidth]{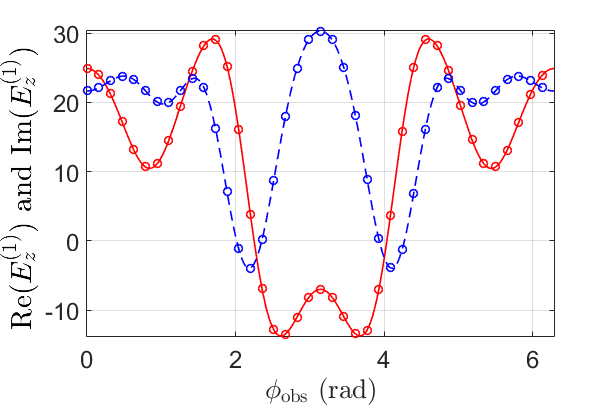}}
    \subfigure[]{\includegraphics[width=0.47\textwidth]{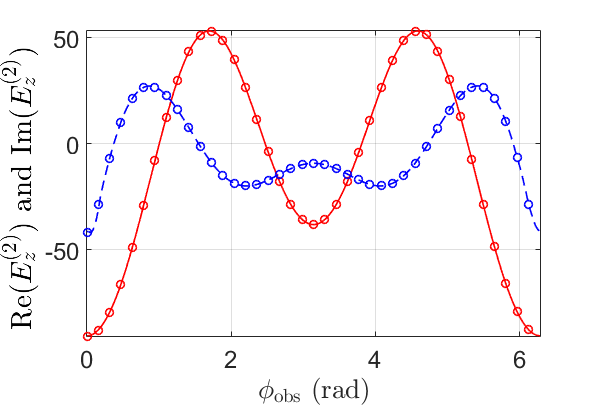}}
    \subfigure[]{\includegraphics[width=0.47\textwidth]{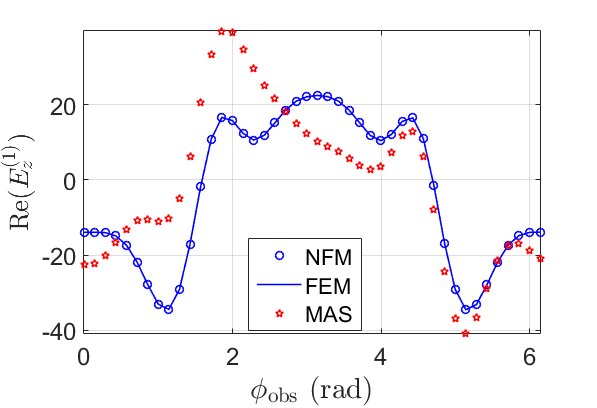}}
    \subfigure[]{\includegraphics[width=0.47\textwidth]{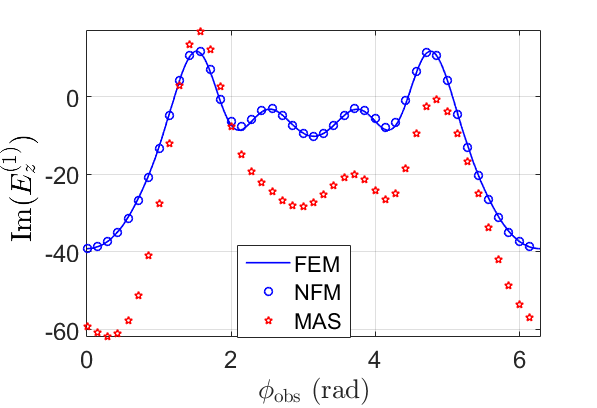}}
    \caption{Real (red solid lines) and imaginary parts (blue dashed lines) of (a) $E_{z,\mathrm{NFM}}^{(1)}(\phi_{\mathrm{obs}})$ and (b) $E_{z,\mathrm{NFM}}^{(2)}(\phi_{\mathrm{obs}})$ for $N=40$, $\sigma^{(1)}_{\mathrm{aux}}=0.33$, $\sigma^{(2)}_{\mathrm{aux}}=5$, $\sigma^{(1)}_{\mathrm{obs}}=5$, and $\sigma^{(2)}_{\mathrm{obs}}=0.5$ along with the corresponding FEM simulations 
    (noted with circles) for external $\rm{TM_z}$ excitation. (c) and (d) are as (a) and (b) but for internal excitation. (e) and (f) illustrate the negative effect of oscillations of the $C^{(1)}_{\mathrm{aux}}$ MAS currents on the computations of the $E^{(1)}_z$ electric field for $N=44$ and external excitation.}
    \label{fig:fields_TM_ell}
\end{figure}

The circular geometry is simple enough to draw analytical and explicit conclusions, which we illustrated numerically in the previous section. 
To demonstrate that those conclusions carry over to noncircular problems, we consider an elliptical boundary (Fig.~\ref{fig:ellipse_EIE}) with  large semi axis $k_1a=2$ and small semi axis $k_1b=1.6$ (i.e., eccentricity $e=0.6$). The excitation source lies on the positive semi-major axis at $k_1 d=4$ for external excitation, and $k_1 d=1$ for internal excitation. The auxiliary surfaces $C^{(1)}_{\mathrm{aux}}$ and $C^{(2)}_{\mathrm{aux}}$ are scaled versions of $C$, with scaling factors $\sigma^{(1)}_{\mathrm{aux}}$ and $\sigma^{(2)}_{\mathrm{aux}}$, respectively. 
Results are shown in Figs.~\ref{fig:dens_TM_ell} and \ref{fig:fields_TM_ell} 
(Fig.~\ref{fig:fields_TM_ell} now compares to FEM simulations 
rather than to exact results). 
It is apparent (compare to Figs.~\ref{fig:dens_vs_norm_TM} and \ref{fig:exact_EIE_el_fields} for the circle) that all results concerning convergence/divergence and existence/nonexistence of oscillations remain unchanged. 
Contrast, in particular,  the correct NFM-SEP fields of Figs.~\ref{fig:fields_TM_ell} (e)-(f) to the corrupted MAS fields in the same figure; this illustrates the negative effect of MAS oscillations (even if these oscillations differ in appearance from those of the circle). 
%
%


\subsection{NFM-SEP vs. MAS: Speed of convergence}

\begin{figure}[H]
    \centering
    \subfigure[]{\includegraphics[width=0.49\textwidth]{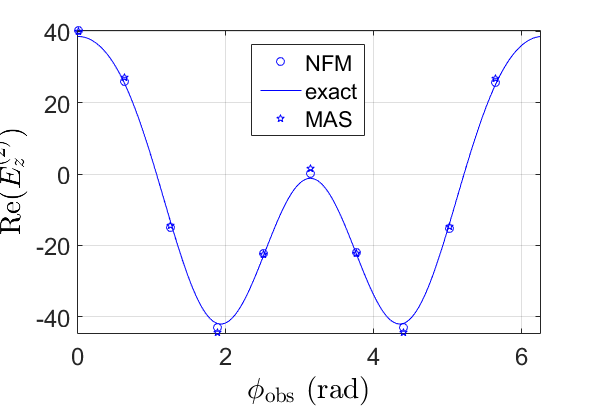}} 
    \subfigure[]{\includegraphics[width=0.49\textwidth]{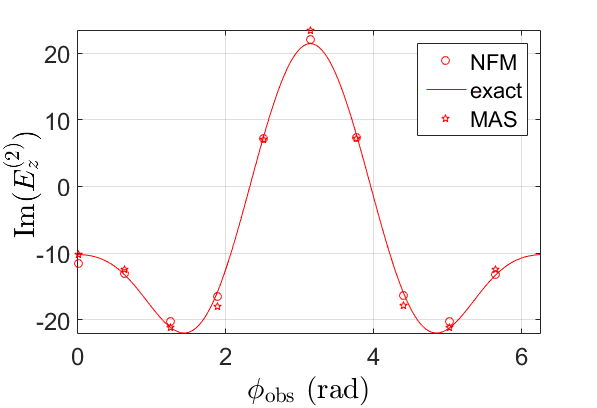}}
    \caption{Real (a) and imaginary (b) parts of  $E^{(2)}_z$ for $\mathrm{TM}_z$ external excitation, as computed from the exact solution (\ref{eq:exact_ext_TM_2}), NFM for $N=10$, and MAS for $N=10$, with $k_1\rho^{(1)}_{\mathrm{aux}}=1$, and $k_1\rho^{(2)}_ {\mathrm{aux}}=4$.}
    \label{fig:int_circ_convergence}
\end{figure}

In view of the comparisons in Sections~\ref{nfm-sep-vs-mas}, \ref{sec:Numerical-circular}, \ref{sec:Numerical-ellipse},
we can ask whether MAS presents any fundamental advantages over NFM-SEP. 
Since the matrix condition numbers of the two methods are identical, the obvious quantity to compare is the speed of convergence. 
Fig.~\ref{fig:int_circ_convergence} shows results from both methods, for $N=10$, together with the exact solution. In this case, NFM-SEP appears slightly better (increasing $N$ to
11 gives MAS results roughly as good as the shown NFM-SEP results). 
More
generally (and as far as the speed of convergence is concerned) any
differences we observed between the two methods were slight (and in some
cases, MAS was slightly better than NFM-SEP). 
A more systematic study of the
convergence speed requires finer measures of comparisons and is
beyond the scope of this paper.
%
%

\section{Conclusions}
\label{sec:conclusions}
A numerical approach based on the combination of the NFM with the SEP was investigated analytically. The combined NFM-SEP was applied for the case of a cylindrical dielectric scatterer excited by an external or internal line source. It was shown that the normalized discrete NFM-SEP currents always converge to the respective continuous surface current densities. This comes in stark contrast to previous findings on the MAS, in which the involved discrete currents diverge and oscillate under certain conditions. To this end, we provided a comparative detailed discussion of NFM-SEP with MAS to highlight many different characteristics of the two methods despite the fact that both methods have the same matrix. The analysis was supplemented by representative numerical results for circular and noncircular shapes verifying the theoretical predictions. An interesting future work direction is to use the NFM-SEP for computing the eigenvalues of dielectric waveguides.

\section*{CRediT authorship contribution statement}

\noindent \textbf{Minas Kouroublakis:} Conceptualization, Methodology, Software, Validation, Formal Analysis, Investigation, Visualization, Writing - Original Draft, Writing - Review \& Editing.\\ \textbf{Nikolaos~L.~Tsitsas:} Conceptualization, Methodology, Investigation, Writing - Review \& Editing.\\ \textbf{George~Fikioris:} Conceptualization, Methodology, Investigation, Writing - Review \& Editing.

\section*{Declaration of Competing Interest}
The authors declare that they have no known competing financial interests or personal relationships that could have appeared to
influence the work reported in this paper.

\section*{Data availability}
No data was used for the research described in the article.

\appendix
\section{Useful Bessel-function formulas}
\label{sec:appendix}
 
This appendix lists certain well-known formulas involving the cylindrical Hankel and Bessel functions \cite{NIST}. Firstly, for $\nu\to+\infty,\, z\in \mathbb{C},\, z\neq 0$, we have the following large-order asymptotic approximations:
\begin{equation}
J_\nu(z)\sim \frac{1}{\sqrt{2\pi \nu}}\left( \frac{e z}{2\nu} \right)^\nu
\label{eq:Bessel},
\end{equation}
\begin{equation}
H^{(2)}_\nu(z)\sim +\mathrm i \sqrt{\frac{2}{\pi \nu}}\left( \frac{e z}{2\nu} \right)^{-\nu}
\label{eq:Hankel},
\end{equation}
\begin{equation}
J'_\nu(z)\sim \sqrt{\frac{\nu}{2\pi}}\left( \frac{e z}{2\nu} \right)^\nu\frac{1}{z}
\label{eq:Bessel_der},
\end{equation}
\begin{equation}
H'^{(2)}_\nu(z)\sim -\mathrm i \sqrt{\frac{2\nu}{\pi}}\left( \frac{e z}{2\nu} \right)^{-\nu}\frac{1}{z}
\label{eq:Hankel_der}.
\end{equation}
Secondly, the addition theorem for the Hankel function $H_0$ is,
\begin{equation}
H_0^{(2)}(\sqrt{x_1^2+x_2^2-2x_1x_2\cos\theta})=\sum_{n=-\infty}^{+\infty}J_n(\min\left\{x_1, x_2  \right\})H_n^{(2)}(\max\left\{ x_1,x_2 \right\})e^{\mathrm i n \theta}
\label{eq:addition_theorem}
\end{equation}
where $x_1>0$, $x_2>0$, and $\theta$ is real. 
For $x_2>x_1$, we obtain the following expressions for the derivatives of (\ref{eq:addition_theorem}) with respect to $x_1$ and $x_2$, respectively,
\begin{align}
\nonumber
&\frac{x_1-x_2\cos \theta}{\sqrt{x_1^2+x_2^2-2x_1x_2\cos\theta}}H_1^{(2)}(\sqrt{x_1^2+x_2^2-2x_1x_2\cos\theta})=\\
&-\sum_{n=-\infty}^{+\infty}J'_n(x_1)H_n^{(2)}(x_2)e^{\mathrm i n \theta},
\label{eq:addition_theorem_der_1}\\
\nonumber
&\frac{x_2-x_1\cos \theta}{\sqrt{x_1^2+x_2^2-2x_1x_2\cos\theta}}H_1^{(2)}(\sqrt{x_1^2+x_2^2-2x_1x_2\cos\theta})=\\
&-\sum_{n=-\infty}^{+\infty}J_n(x_1)H'^{(2)}_n(x_2)e^{\mathrm i n \theta}.
\label{eq:addition_theorem_der_2}
\end{align}
Finally, we note the following Wronskian identity,
\begin{equation}
\label{eq:wronskian}
    J_n(z)H_n^{'(2)}(z)-J'_n(z)H_n^{(2)}(z)=\frac{2}{\mathrm i \pi z}.
\end{equation}
%






\section{Verification of (\ref{eq:current_densities_ext_1}) and (\ref{eq:current_densities_ext_2})}
\label{sec:correct-fields}

In this appendix, we verify that the surface current densities in 
(\ref{eq:current_densities_ext_1}) and (\ref{eq:current_densities_ext_2})  yield the correct fields in (\ref{eq:exact_ext_TM_1}) and (\ref{eq:exact_ext_TM_2}). For the field in $R_1$, use the addition theorems (\ref{eq:addition_theorem})-(\ref{eq:addition_theorem_der_2}) to rewrite  (\ref{eq:R1_field})  as follows:
\begin{multline}
E_z^{(1)}(\rho_{\mathrm{obs}}, \phi_{\mathrm{obs}})=-\frac{k_1Z_1}{4}IH_0^{(2)}(k_1D_{\mathrm{fil,obs}})\\
-\frac{k_1Z_1\rho_{\mathrm {cyl}}}{4}\sum_{n=-\infty}^{+\infty}J_n(k_1\rho_{\mathrm{cyl}})H_n^{(2)}(k_1\rho_{\mathrm {obs}})e^{\mathrm i n \phi_{\mathrm {obs}}}\int_{-\pi}^{\pi}J^s_z(\phi_{\mathrm{cyl}})e^{-\mathrm in \phi_{\mathrm{cyl}}}\mathrm{d}\phi_{\mathrm{cyl}}\\-\frac{k_1\rho_{\mathrm{cyl}}}{4\mathrm i}\sum_{n=-\infty}^{+\infty}J'_n(k_1\rho_{\mathrm{cyl}})H_n^{(2)}(k_1\rho_{\mathrm {obs}})e^{\mathrm i n \phi_{\mathrm {obs}}}\int_{-\pi}^{\pi}M^s_\phi(\phi_{\mathrm{cyl}})e^{-\mathrm in \phi_{\mathrm{cyl}}}\mathrm{d}\phi_{\mathrm{cyl}}
\label{eq:R1_field_2}
\end{multline}
The integrals are $2\pi$ times the $n$-th Fourier series
coefficient $I^s_n$ (of $J_z^s$) and $M^s_n$ (of $M_\phi^s$). Thus, by (\ref{eq:current_densities_ext_1}) and (\ref{eq:current_densities_ext_2}), we have
\begin{multline}
E_z^{(1)}(\rho_{\mathrm{obs}}, \phi_{\mathrm{obs}})=-\frac{k_1Z_1}{4}IH_0^{(2)}(k_1D_{\mathrm{fil,obs}})\\-\frac{k_1Z_1I}{4}\sum_{n=-\infty}^{+\infty}J_n(k_1\rho_{\mathrm{cyl}})H_n^{(2)}(k_1\rho_{\mathrm {obs}})\times \\\frac{Z_1H_n^{(2)}(k_1\rho_{\mathrm{fil}})J'_n(k_2\rho_{\mathrm{cyl}})}{Z_1H_n^{(2)}(k_1\rho_{\mathrm{cyl}})J'_n(k_2\rho_{\mathrm{cyl}})-Z_2H_n^{'(2)}(k_1\rho_{\mathrm{cyl}})J_n(k_2\rho_{\mathrm{cyl}})}e^{\mathrm i n \phi_{\mathrm {obs}}}\\-\frac{k_1I}{4\mathrm i}\sum_{n=-\infty}^{+\infty}J'_n(k_1\rho_{\mathrm{cyl}})H_n^{(2)}(k_1\rho_{\mathrm {obs}})\times\\\frac{\mathrm i Z_1Z_2H_n^{(2)}(k_1\rho_{\mathrm{fil}})J_n(k_2\rho_{\mathrm{cyl}})}{Z_1H_n^{(2)}(k_1\rho_{\mathrm{cyl}})J'_n(k_2\rho_{\mathrm{cyl}})-Z_2H_n^{'(2)}(k_1\rho_{\mathrm{cyl}})J_n(k_2\rho_{\mathrm{cyl}})}e^{\mathrm in \phi_{\mathrm{obs}}},
\end{multline}
which after some simple manipulations yields the exact solution (\ref{eq:exact_ext_TM_1}). Application of the same procedure to (\ref{eq:R2_field}) yields 
\begin{multline}
 E_z^{(2)}(\rho_{\mathrm{obs}}, \phi_{\mathrm{obs}})= -\frac{k_2Z_2I}{4}\sum_{n=-\infty}^{+\infty}J_n(k_2\rho_{\mathrm{obs}})H_n^{(2)}(k_2\rho_{\mathrm{cyl}})\times \\\frac{Z_1H_n^{(2)}(k_1\rho_{\mathrm{fil}})J'_n(k_2\rho_{\mathrm{cyl}})}{Z_1H_n^{(2)}(k_1\rho_{\mathrm{cyl}})J'_n(k_2\rho_{\mathrm{cyl}})-Z_2H_n^{'(2)}(k_1\rho_{\mathrm{cyl}})J_n(k_2\rho_{\mathrm{cyl}})}e^{\mathrm i n \phi_{\mathrm {obs}}}\\-\frac{k_2I}{4\mathrm i}\sum_{n=-\infty}^{+\infty}J_n(k_2\rho_{\mathrm {obs}})H_n^{'(2)}(k_2\rho_{\mathrm{cyl}})\times \\\frac{\mathrm i Z_1Z_2H_n^{(2)}(k_1\rho_{\mathrm{fil}})J_n(k_2\rho_{\mathrm{cyl}})}{Z_1H_n^{(2)}(k_1\rho_{\mathrm{cyl}})J'_n(k_2\rho_{\mathrm{cyl}})-Z_2H_n^{'(2)}(k_1\rho_{\mathrm{cyl}})J_n(k_2\rho_{\mathrm{cyl}})}e^{\mathrm i n \phi_{\mathrm {obs}}},
\end{multline}
which after some manipulations and use of the identity (\ref{eq:wronskian}), yields (\ref{eq:exact_ext_TM_2}).

\bibliographystyle{elsarticle-num}
\bibliography{mybib_rev_kats}

\end{document}